\documentclass{amsart}
\usepackage{amsmath,amssymb,amsfonts}
\usepackage{setspace}
\usepackage[dvipsnames]{xcolor}

\newcommand*{\mailto}[1]{\href{mailto:#1}{\nolinkurl{#1}}}

\usepackage{amssymb}
\usepackage{hyperref}
\hypersetup{
  colorlinks,
  linkcolor={red!50!black},
  citecolor={blue!50!black},
  urlcolor={blue!80!black}
}
\usepackage{graphicx}
\usepackage{curves}
\usepackage{enumerate}
\usepackage{tikz}
\usetikzlibrary{matrix}
\usepackage{marginnote}

\let\phi=\varphi
\let\circc=\circ
\renewcommand{\circ}{\hspace*{-0.5mm}\circc\hspace*{-0.5mm}}
\newcommand{\diff}{\text{diff}(\Real)}
\newcommand{\brho}{{\bar\rho}}
\newcommand{\br}{{\bar r}}
\newcommand{\bs}{{\bar s}}

\newcommand{\fracpar}[2]{\frac{\partial #1}{\partial #2}}
\newcommand{\fracdelta}[2]{\frac{\delta #1}{\delta #2}}

\newcommand{\epsi}{\varepsilon}

\newcommand{\Helm}{\mathcal{H}}
\newcommand{\invHelm}{\Helm^{-1}}

\newcommand{\Ekin}{\ensuremath{E^{\text{kin}}} }
\newcommand{\Epot}{\ensuremath{E^{\text{pot}}} }

\newcommand{\Lcal}{\mathcal{L}}
\newcommand{\LcalCH}{\Lcal_{\text{CH}}}

\newcommand{\Lpot}{\ensuremath{\mathcal{L}^{\text{pot}}} }

\newcommand{\Etot}{\ensuremath{E^{\text{tot}}} }
\newcommand{\Real}{\mathbb R}

\newcommand{\abs}[1]{\left|#1\right|}

\DeclareMathOperator{\id}{id}

\newcommand{\brac}[1]{\left<#1\right>}
\DeclareMathOperator{\sign}{sign}

\newtheorem{theorem}{Theorem}[section]

\newtheorem{definition}[theorem]{Definition}

\numberwithin{equation}{section}
\usepackage{soul}
\soulregister\cite7
\soulregister\ref7
\soulregister\pageref7
\soulregister\eqref7
\soulregister\xaviername7
\soulregister\katrinname7
\sethlcolor{yellow}

\def\xaviername{\textcolor{ForestGreen}{\bf Xavier}}
\def\katrinname{\textcolor{blue}{\bf Katrin}}

\begin{document}

\title[Symmetries and multipeakons for the M2CH system]{Symmetries and
  multipeakon solutions for the modified two-component Camassa--Holm system}
  
 \author[K. Grunert]{Katrin Grunert}
\address{Department of Mathematical Sciences\\NTNU\\
  Norwegian University of Science and Technology\\
  7491 Trondheim\\ Norway} \email{\mailto{katrin.grunert@ntnu.no}}

\author[X. Raynaud]{Xavier Raynaud}
\address{Applied Mathematics, SINTEF ICT, Oslo, Norway\\
  and Department of Mathematical Sciences\\NTNU\\
  Norwegian University of Science and Technology\\
  7491 Trondheim, Norway} \email{\mailto{xavier.raynaud@ntnu.no}}

\thanks{Research supported by the grant {\it Waves and Nonlinear Phenomena (WaNP)} from the Research Council of Norway.}

\dedicatory{This paper is dedicated to Helge Holden on the occasion of his sixtieth anniversary, with admiration and gratefulness for all the inspiration
  he has been giving us in our work.}

\begin{abstract}
  Compared with the two-component Camassa--Holm system, the modified two-component
  Camassa--Holm system introduces a regularized density which makes possible the
  existence of solutions of lower regularity, and in particular of multipeakon
  solutions. In this paper, we derive a new pointwise invariant for the modified
  two-component Camassa--Holm system. The derivation of the invariant uses
  directly the symmetry of the system, following the classical argument of
  Noether's theorem. The existence of the multipeakon solutions can be directly
  inferred from this pointwise invariant. This derivation shows the strong
  connection between symmetries and the existence of special solutions. The
  observation also holds for the scalar Camassa--Holm equation and, for
  comparison, we have also included the corresponding derivation. Finally, we
  compute explicitly the solutions obtained for the peakon-antipeakon case. We
  observe the existence of a periodic solution which has not been reported in
  the literature previously. This case shows the attractive effect that the
  introduction of an elastic potential can have on the solutions.
\end{abstract}

\maketitle

\section{Introduction}

In \cite{holm2009singular}, the authors introduce the modified two-component
Camassa--Holm system (M2CH), which is given by
\begin{subequations}
\label{eq:m2chfull}
  \begin{align}
    \label{eq:m2chfull1}
    m_t + u m_x + 2 m u_x + \brho_x\rho  &=0,\\
    \label{eq:m2chfull2}
    \rho_t + (u \rho)_x &=0,
  \end{align}
  where
  \begin{align}
    \label{eq:m2chfull3}
    m &= u - u_{xx},\\
    \label{eq:m2chfull4}
    \rho &= \brho - {\brho}_{xx}.
  \end{align}
\end{subequations}
This system originates from the Camassa--Holm (CH) equation,
\begin{equation}
  \label{eq:ch}
  m_t + u m_x + 2 m u_x = 0,
\end{equation}
with \eqref{eq:m2chfull3} and the two-component Camassa-Holm system (2CH)
\begin{equation}
  \label{eq:chsys}
  m_t + u m_x + 2 m u_x + \rho_x\rho = 0,
\end{equation}
with \eqref{eq:m2chfull2} and \eqref{eq:m2chfull3}. All these equations can be
derived from a variational principle for the kinetic energy that is defined
\begin{equation}
  \Ekin(t) = \frac12\int_\Real (u^2+u_x^2)(t,x)\,dx,
\end{equation}
and the following potential energy
\begin{equation}
  \label{eq:defallEpot}  
  \Epot = 0,\quad \Epot = \frac12\int_\Real\rho^2(t,x)\,dx,\quad \Epot = \frac{1}{2}\int_\Real (\brho^2 + \brho_x^2)(t,x)\,dx,
\end{equation}
for CH, 2CH and M2CH, respectively. An advantage of M2CH is that the system of
equations requires a lower regularity for the density, compared to 2CH. Indeed,
given the potential energy as in \eqref{eq:defallEpot}, while the 2CH
system requires that $\rho\in L^2(\Real)$, the M2CH system requires that
$\brho\in H^1(\Real)$, which is equivalent to $\rho\in H^{-1}(\Real)$, as the
Helmholtz operator $\id - \partial_{xx}$ is an isomorphism from $H^1(\Real)$ to
$H^{-1}(\Real)$.

The CH equation has a rich mathematical structure which explains the very
extensive literature that is available on this equation. In this work, we
consider global conservative solutions which can be defined beyond the blow-up
of the classical solutions. For the CH equation, the blow-up scenario is known
and occurs when, for some given initial data $u_0\in H^1(\Real)$, the spatial
derivative $u_x$ becomes unbounded from below within finite time, while the
$H^1(\Real)$-norm of $u$ and hence also $u$ remain bounded. This phenomenon,
which is referred to as wave breaking is described in
\cite{constantin2000,constantinescher1998a,
  constantinescher1998b,constantinescher2000}. In particular, it
can be predicted if wave breaking occurs in the nearby future or not, see
\cite{grunert2015break} and the references therein.  In more recent works, the
regularization properties of the density in the case of the 2CH system have been
studied \cite{constantin2008integrable,grunert2012global,guan2010global}. There,
it is shown that if the density is bounded away from zero initially, a solution
with smooth initial data will never experience blow-up. We find appropriate the
following interpretation. The governing equations, that are obtained from the
variational principle, model the velocity $u$ of an underlying flow map
$\phi(t,\xi)$, that is $\phi_t = u\circ\phi$. The elastic energy introduced by
$\Epot$ prevents compression so that the flow map cannot become irregular in the
sense that several particles can occupy the same place $\phi(t,\xi_1) =
\phi(t,\xi_2)$ for two \textit{particles} $\xi_1$ and $\xi_2$. The potential
energy for the M2CH system is weaker in the sense that, if we consider a
concentration of particles at a single point, the potential energy for the 2CH
system becomes infinite making this state not reachable while it is finite for
the M2CH system. Indeed, formally speaking, a concentration of particles gives
rise to a density $\rho$ equal to a Dirac delta function which has infinite
$L^2(\Real)$ norm while $\brho=\frac12e^{-\abs{x}}\star \rho$ remains in
$H^1(\Real)$. However, compared to the 2CH system, the M2CH system has the
property of having a special class of solutions.

\begin{figure}[h]
  \centering
  \includegraphics[width=\textwidth]{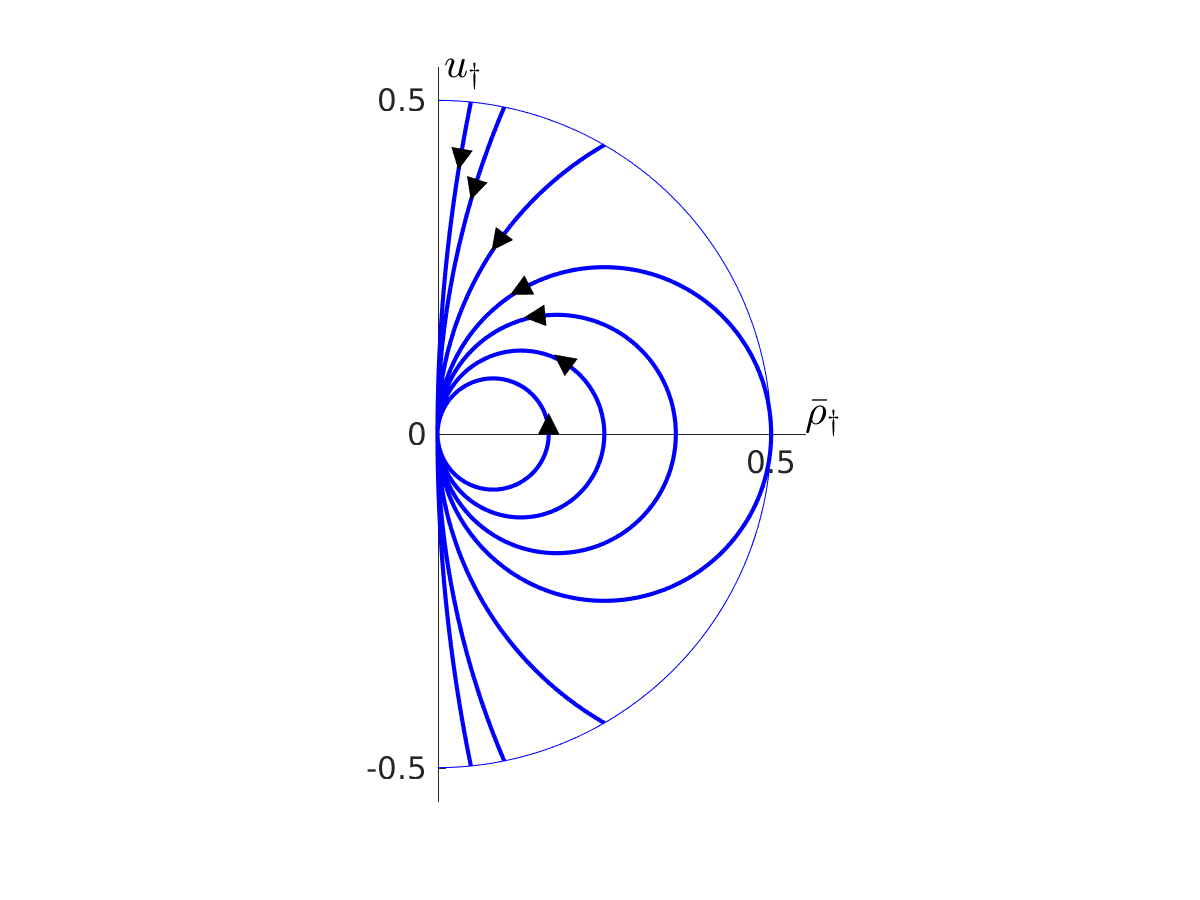}
  \caption{Plot of the trajectory of $(\brho_\dagger, u_\dagger)$ for different
    values of $s$. The outer half-circle represents the limiting values of the
    solution when $t\to\pm\infty$ when $s\leq 1$, see \eqref{eq:limcirc}. The
    circles in the middle represent the periodic solution for $s\geq 1$.}
  \label{fig:trajplot}
\end{figure}

The CH equation admits a special type of soliton-like solutions that have been
called multipeakons, due to the peaks that characterize them. The multipeakon
solutions can be seen as a discrete version of the equation. Such solutions are
dense \cite{BressanFonte2005}, robust \cite{constantinstrauss2000, dikamolinet2009, dikamolinet20092}, and have been used to design convergent numerical
schemes, which can also handle blow-up
\cite{HoldenRaynaud2006,HoldenRaynaud2006b}. It turns out that the M2CH also
admits such solutions, as pointed out in \cite{holm2009singular}. In this paper,
we follow the following understanding. Special solutions exist because the
equation has a special structure and structures are identified by symmetries. In
this case, the symmetry of the system is related to the invariance with respect
to relabeling of both the kinetic and potential energy. From Noether's theorem,
we know that this invariance must imply the existence of conservation
laws. Since the group of diffeomorphism has infinite dimension, we expect
infinitely many invariants. As we will see, the Noether argument leads us to
pointwise invariants of the form
\begin{equation}
  \label{eq:invpres}
  (u - u_{xx})(t,\phi(t,\xi))\phi_{\xi}^2(t,\xi)
\end{equation}
in the case of the CH equation, which also encodes the conservation of the left angular momentum \cite{kolev2004}. Such invariants have been derived much earlier
in \cite{arnold1999topological} but here, we present a more straightforward
derivation that does not require the advanced topological framework used in the
fore-mentioned work. Of course, we miss some fundamental insight but simplifying
the derivation, we can make it possible to adapt it directly to the case of the
M2CH system later. The problem of the pointwise invariant \eqref{eq:invpres} is
that it is not so easy to exploit as it mixes natural Eulerian variables (the
expression of $u_{xx}$ is complicated in Lagrangian variables) and Lagrangian
variables ($\phi$ it not directly available from Eulerian variables). However,
it can be used to show the existence of multipeakon solutions, thus making clear
the connection between the symmetries of the system and the existence of a large
and non-trivial class of special solutions. In this paper, we modify the
variational formulation of the M2CH system to make it suitable for the use of
the Noether's argument. We derive the pointwise invariant of the M2CH system and
describe how the existence of multipeakon solutions can be inferred from it. The
variational formulations are always done with respect to the flow map. Hence,
Lagrangian variables are naturally introduced in this setting. The change of
variables to Lagrangian variables is known to be a mean of getting rid of
non-linearity in the advection term corresponding to $u_t + uu_x$ in the
equation below. We denote by $\Helm=\partial - \partial_{xx}$ the Helmholtz
operator. After applying its inverse $\invHelm$ to \eqref{eq:m2chfull}, this
system of equations becomes
\begin{subequations}
  \label{eq:m2ch}
  \begin{align}
    \label{eq:m2ch1}
    u_t + uu_x + P_x &= 0\\
    \label{eq:m2ch2}
    \brho_t+u\brho_x+R+S_x&=0
\end{align}
with
\begin{align}
    \label{eq:m2ch3}
    P  &= \invHelm \left(u^2 + \frac12u_x^2 + \frac12\brho^2 - \frac12\brho_x^2\right)\\
    \label{eq:m2ch4}
    R& = \invHelm\left(u_x\brho\right) \\
    \label{eq:m2ch5}
    S&= \invHelm\left(u_x\brho_x\right)
  \end{align}
\end{subequations}
The Lagrangian variables are given by the characteristics defined as $y_t(t,\xi)
= u(t,y(t,\xi))$ and the Lagrangian velocity $U(t,\xi) = u(t,y(t,\xi))$. In the
case of the M2CH system, we need to introduce more variables to be able to
handle the blow-up of the solution.  Here, we follow the approach presented in
\cite{guantanwei2015,tan2011global} which is very close the one introduced in
\cite{grunert2012global,holdenraynaud2007}. Once the system of equations is
completely rewritten in term of purely Lagrangian variables, semi-linearities in
the system enable us to obtain global solutions. Thus, the Lagrangian system
defines the solutions, which are then mapped back to Eulerian variables in order
to obtain some weak solutions to the original M2CH system. We show that the
existence of the pointwise invariant implies the existence of multipeakon solutions
but, even if this invariant can be expressed in term of purely Lagrangian
variables, its form becomes then intricate. However, this fundamental invariance
property is preserved by the change of Lagrangian variables so that the
existence of multipeakon can also be obtained in the Lagrangian setting, see
\cite{holdenraynaud2007b} for the corresponding work in the case of the CH
equation. In section \ref{sec:doublemultipeak}, we derive the multipeakon
equations directly from the system \eqref{eq:m2ch} and not as in
\cite{holm2009singular} where a discretization of the Hamiltionian is used. We
compute explicit solutions in the case of the anti-symmetric peakon-antipeakon
solution. We discover an interesting dynamic in this case, which can be
decomposed in three different cases. For all cases, the peaks collide but there
are different behaviors when the peaks move away from each other, after
collision. In the first case, there is not enough potential energy in the system
to retain the particles from completely departing from each other. In the second
case, the potential energy prevents them to do so. We can compare the situation
to a classical discrete mechanical system where Hooke's law is used to model the
elastic forces. Such elastic forces act in both ways. They are repulsive when
particles approach each other, over a given equilibrium state, and attractive
when the particle move far away from each other. For M2CH, we observe that the
potential energy does not yield a repulsive force that is strong enough to
prevent collision but its attractive effect can prevent the fully departure of
the peaks from each others. The solution in this case is periodic and we finally
end up with a oscillatory system where the kinetic energy and the potential
energy vanishes one after the other, as for a standard pendulum. This type of
solution has not been observed for CH or 2CH. The last case in the description
of the dynamics is the limiting case, where we do not obtain a periodic solution
but the peaks are slowed down by the attractive force until their velocity
vanishes. The position of the (left) right particle tends to (minus) infinity
while their velocity tend to zero, see Figure \ref{fig:trajplot}.
  
\section{Conservation laws}
\label{sec:conslaws}

For the M2CH system, we define the \textit{kinetic energy} as
\begin{equation}
  \label{eq:defEkin}
  \Ekin(t) = \frac12\int_\Real (u^2+u_x^2)(t,x)\,dx.
\end{equation}
The proper definition of kinetic energy from physics is $\frac12\int_\Real\rho
u^2\,dx$. However, we are going to see that the term defined in
\eqref{eq:defEkin} plays a role which resembles the one of the kinetic energy in
standard physical systems and that is why we use this terminology. Using the same
type of analogy, we refer to the quantity $\Epot$ defined below as the
\textit{potential energy},
\begin{equation}
  \label{eq:defEpot}
  \Epot(t) = \frac{1}{2}\int_\Real (\brho^2 + \brho_x^2)(t,x)\,dx.
\end{equation}
The total energy is then given by
\begin{equation}
  \label{eq:defEtot}
  \Etot(t) = \frac12\int_\Real (u^2 + u_x^2 + \brho^2 + \brho_x^2)(t,x)dx.
\end{equation}
The M2CH system can be derived from a variational principle
for the Lagrangian
\begin{equation}
  \label{eq:defLm2ch}
  \Lcal = \Ekin - \Epot.
\end{equation}
We do not give the details for this computation and refer instead to
\cite{grunert2016}. The invariance of the Lagrangian with respect to time
implies through Noether's theorem that the total energy as defined in
\eqref{eq:defEtot} is preserved in time. More precisely, we have the following
conservation law for the energy,
\begin{multline}
  \label{eq:energconslaw}
  (u^2 + u_x^2 + \brho^2 + \brho_x^2)_t + (u(u^2 + u_x^2 + \brho^2 + \brho_x^2))_x = \\
  \left(u^3 - 2Pu - 2\brho\invHelm(\rho u) + 2u\bar\rho^2\right)_x.
\end{multline}
Let us derive \eqref{eq:energconslaw} from \eqref{eq:m2ch}. One can prove that
for any smooth function $q(t,x)$ which satisfies
\begin{equation}
  \label{eq:q}
  q_t + uq_x + Q = 0
\end{equation}
for some given smooth $Q$, one has
\begin{equation}
  \label{eq:qqx}
  (q^2 + q_x^2)_t + (u(q^2 + q_x^2))_x = -q_x^2u_x + u_xq^2 - 2qQ - 2q_xQ_{x}.
\end{equation}
We let the reader check this property. From \eqref{eq:m2ch2}, we have
\begin{equation*}
  \brho_t + u\brho_x  + \invHelm(\rho u )_x - u\invHelm\rho_x = 0.
\end{equation*}
Let us define
\begin{equation*}
  Q = \invHelm(\rho u )_x - u\invHelm\rho_x.
\end{equation*}
After some computations, we get
\begin{equation*}
  Q_x = \invHelm(\rho u) - u_x\brho_x - u\brho.
\end{equation*}
Hence,
\begin{equation*}
  \brho_x Q_x + \brho Q = (\brho\invHelm(\rho u))_x - u_x\brho_x^2 - 2u\brho\brho_x.
\end{equation*}
Then, we apply \eqref{eq:qqx} for $q = \brho$ and get
\begin{equation}
  \label{eq:brhoconslaw}
  (\brho^2 + \brho_x^2)_t + (u(\brho^2 + \brho_x^2))_x = -2(\brho\invHelm(\rho u))_x + u_x(\brho_x^2 - \brho^2) + 2(\brho^2 u)_x
\end{equation}
Now, we set
\begin{equation*}
  Q = P_x
\end{equation*}
and we apply \eqref{eq:qqx} for $q = u$ and get
\begin{equation}
  \label{eq:uconslaw}
  (u^2 + u_x^2)_t + (u(u^2 + u_x^2))_x = (u^3 - 2uP)_x + 2u_x(\Helm P - u^2 - \frac12u_x^2)
\end{equation}
We sum up \eqref{eq:brhoconslaw} and \eqref{eq:uconslaw} and obtain
\begin{multline*}
  (u^2 + u_x^2 + \brho^2 + \brho_x^2)_t + (u(u^2 + u_x^2 + \brho^2 + \brho_x^2))_x = \\
  \left(u^3 - 2Pu - 2\brho\invHelm(\rho u) + 2u\bar\rho^2\right)_x +
  2u_x\left(\Helm P - u^2 - \frac12u_x^2 + \frac12\brho_x^2 - \frac12\brho^2\right).
\end{multline*}
The last of the two terms on the right-hand side vanishes because of
\eqref{eq:m2ch3} so that the conservation law \eqref{eq:energconslaw} follows.
  
\section{Lagrangian variables}
\label{sec:lagvar}
In this section, we describe how the M2CH system \eqref{eq:m2ch} can be
rewritten in Lagrangian variables to obtain a system which is formally
equivalent but whose linear structure can be used to prove the global existence
of solutions. In this section the derivation of the equivalent system is only
formal. Once the system is obtained, the construction of the solution in
Lagrangian variables and the mapping back to the original Eulerian variables can
be done rigorously, see \cite{guantanwei2015}. We introduce the characteristics
defined as
\begin{equation}
  \label{eq:defy}
  y_t(t,\xi) = u(t,y(t,\xi)),
\end{equation}
the \textit{Lagrangian velocity} defined as
\begin{equation}
  \label{eq:defU}
  U(t,\xi) = u(t, y(t,\xi)),
\end{equation}
the cumulative total \textit{energy distribution} defined as
\begin{equation}
  \label{eq:defH}
  H(t,\xi) = \int_{-\infty}^{y(t,\xi)}(u^2 + u_x^2 + \brho^2 + \brho_x^2)(t,x)\,dx,
\end{equation}
the \textit{Lagrangian regularized potential energy} defined as
\begin{equation}
  \label{eq:defbr}
  \br(t, \xi) = \brho(t, y(t,\xi)),
\end{equation}
and, finally,
\begin{equation}
  \label{eq:defbs}
  \bs (t,\xi) = \brho_x(t, y(t,\xi)).
\end{equation}
We will assume in the remaining of this formal derivation that the derivative
$y_\xi$ does not vanish. After differentiating \eqref{eq:defbr} and
\eqref{eq:defbs}, we observe that, formally,
\begin{equation}
  \label{eq:relbrxibs}
  \br_\xi=\bs y_\xi
\end{equation}
and
\begin{equation}
  \label{eq:relrbrho}
  \bar s_\xi = \br y_\xi - \rho(t, y)y_\xi.
\end{equation}
The inverse Helmholtz operator can be written using Green's function as
\begin{equation}
  \label{eq:invHelmExpl}
  [\invHelm q](x) = \frac12\int_\Real e^{-\abs{x -z}}q(z)\,dz.
\end{equation}
Hence, we get from the definition \eqref{eq:m2ch3} of $P$ that
\begin{equation*}
  P(t,x) = \frac14\int_\Real e^{-\abs{x-z}}((u^2 + u_x^2 + \brho^2 + \brho_x^2) + (u^2 -2 \brho_x^2))(t,z)\,dz.
\end{equation*}
We change to Lagrangian variables and use the definition \eqref{eq:defH} of $H$ and the identity \eqref{eq:defbs} to get
\begin{equation}
  \label{eq:defLP}
  P(t,y) = \frac14\int_{\Real} e^{-\abs{y(t,\xi) - y(t,\eta)}}(H_\xi(t,\eta) + (U^2(t,\eta) - 2\bs^2(t,\eta))y_\xi(t,\eta))\,d\eta.
\end{equation}
The change to Lagrangian variables has the recognized advantage to get rid of the
first non-linear term in \eqref{eq:m2ch1}, which becomes
\begin{equation}
  \label{eq:UtPx}
  U_t(t,\xi) = -P_x(t, y(t,\xi))
\end{equation}
Introducing
\begin{multline}
  \label{eq:defQ}
  Q(t,\xi) = -\frac14\int_\Real\sign(\xi - \eta) e^{-\abs{y(t,\xi) - y(t,\eta)}}\Big(H_\xi(t,\eta) \\
  + (U^2(t,\eta) - 2\bs^2(t,\eta))y_\xi(t,\eta)\Big)\,d\eta,
\end{multline}
and assuming that $y_\xi$ remains strictly positive, we can differentiate $P$ in \eqref{eq:defLP} and obtain that
\begin{equation*}
  P_x(t, y(t, \xi)) y_\xi(t, \xi) = Q(t, \xi)y_{\xi}(t, \xi).
\end{equation*}
We simplify the above expression by $y_\xi$ and thus \eqref{eq:UtPx} yields
\begin{equation}
  \label{eq:laggovUt}
  U_t = - Q.
\end{equation}
Following the same lines we introduce the integrated variables
\begin{subequations}\label{eq:defRSVW}
\begin{align}
R(t,\xi) & =\frac12 \int_\Real e^{-\abs{y(t,\xi)-y(t,\eta)}}U_\xi(t,\eta) \br(t,\eta)\,d\eta,\\
S(t,\xi)& = \frac12 \int_\Real e^{-\abs{y(t,\xi)-y(t,\eta)}} U_\xi(t,\eta)\bs(t,\eta)\,d\eta,\\
V(t,\xi) & =-\frac12 \int_\Real \sign(\xi-\eta)e^{-\abs{y(t,\xi)-y(t,\eta)}}U_\xi(t,\eta) \br(t,\eta)\,d\eta,\\
W(t,\xi)& = -\frac12 \int_\Real \sign(\xi-\eta)e^{-\abs{y(t,\xi)-y(t,\eta)}} U_\xi(t,\eta)\bs(t,\eta)\,d\eta,
\end{align}
\end{subequations}
and we obtain that
\begin{equation}
  \label{eq:derRS}
  R_\xi = Vy_\xi \quad\text{ and }\quad  S_\xi = Wy_\xi.
\end{equation}
Moreover, after differentiation, we get
\begin{equation}
  \label{eq:derVW}
  V_\xi = -U_\xi\br + Ry_\xi \quad\text{ and }\quad W_\xi = -U_\xi\bs + Sy_\xi.
\end{equation}
The conservation law \eqref{eq:energconslaw} gives 
\begin{equation}\label{help:help}
H_t= \left(u^3-2Pu-2\bar\rho\invHelm(\rho u)+2u\brho^2\right)(t,y).
\end{equation}
Recalling that $\rho=\brho-\brho_{xx}$, direct computations yield that
\begin{equation*}
\brho\invHelm(\rho u)= \brho\invHelm(\brho u-(\brho u)_{xx}+(u_x\brho)_x+u_x\brho_x)=u\brho^2+\brho R_x+\brho S
\end{equation*}
and \eqref{help:help} can be rewritten as 
\begin{equation}
  \label{eq:govHLag}
  H_t = U^3 - 2PU - 2\br (S+V).
\end{equation}
Note that in \eqref{eq:govHLag}, we slightly abused the notations and denoted
$P(t,y)$ as $P(t,\xi)$. We continue to do so in the remaining.  Differentiating
$P$ and $Q$ gives us
\begin{align}
  \label{eq:Pxi}
  P_\xi &= Q y_\xi,\\
  \label{eq:Qxi}
  Q_\xi &= -\frac12H_\xi - (\frac12 U^2-\bs^2 - P)y_\xi.
\end{align}
For the Lagrangian regularized potential energy density, \eqref{eq:m2ch2} yields
\begin{equation}
  \label{eq:govrlag}
  \br_t = - (R + W).
\end{equation}
Let us now gather the governing equations we have obtained in \eqref{eq:defy},
\eqref{eq:laggovUt}, \eqref{eq:govHLag}, and \eqref{eq:govrlag}. We have seen
that the governing equations \eqref{eq:m2ch} are formally equivalent to the
system
\begin{subequations}
    \label{eq:laggoveq}
  \begin{align}
    \label{eq:laggoveq1}
    y_t &= U,\\
    \label{eq:laggoveq2}
    U_t &= -Q,\\
    \label{eq:laggoveq3}
    H_t &= U^3 - 2PU - 2\br (S+V),\\
    \label{eq:laggoveq4}
    \br_t &= - (R + W),\\
    \label{eq:laggoveq5}
    \bs_t& = - (S + V),
  \end{align}
\end{subequations}
where the quantities $P$, $Q$, $R$, $S$, $V$, and $W$ are defined in
\eqref{eq:defLP}, \eqref{eq:defQ}, and \eqref{eq:defRSVW}, respectively. We can
differentiate the first four equations in \eqref{eq:laggoveq} and obtain
\begin{subequations}
    \label{eq:laggoveqder}
  \begin{align}
    \label{eq:laggoveqder1}
    y_{\xi,t} &= U_\xi,\\
    \label{eq:laggoveqder2}
    U_{\xi,t} &= \frac12H_\xi + \Big(\frac12 U^2-\bs^2 - P\Big)y_\xi,\\
    \label{eq:laggoveqder3}
    H_{\xi,t} &= \Big(3U^2 - 2P+2\br^2 \Big) U_\xi - 2\Big(QU +\br(V + W)\Big)y_\xi - 2\Big(S+V\Big)\br_\xi,\\
    \label{eq:laggoveqder4}
    \br_{\xi,t} &= \bs U_\xi -(S+V) y_\xi.
  \end{align}
\end{subequations}
The system \eqref{eq:laggoveqder} reveals the semi-linear nature of the
equivalent system. Indeed, the system is semi-linear with respect to the
derivatives $y_\xi$, $U_\xi$, $H_\xi$ and $\br_\xi$ in the sense that all the
other terms (included $\bs$) that enter the system are of higher regularity than
these derivatives. The semi-linearity of the system is essential in the proof of
the existence of solutions using Picard's argument. The variable $H$ is now
considered as an independent variable but when we introduced it in
\eqref{eq:defH}, it was clearly dependent on the other variables. Changing
variables in \eqref{eq:defH} gives us
\begin{equation}
  \label{eq:enidlag}
  y_\xi H_\xi = (U^2 + \br^2 +\bar s^2)y_\xi^ 2 + U_\xi^2,
\end{equation}
and it can be shown that the governing system \eqref{eq:laggoveq} preserves this
identity, if it holds initially. Thus, we have decoupled $H$ from the other
variables, in particular to obtain a semi-linear system, but \eqref{eq:enidlag}
shows that the variables are not truly independent as they are constrained by
the system to remain on the ``manifold'' defined by \eqref{eq:enidlag}. As shown
in \cite{guantanwei2015}, the system of ordinary differential equations
\eqref{eq:laggoveq} has global solutions in a suitable Banach space. In
particular,
\begin{equation}\label{space:solution}
  y-\id, H\in L^\infty(\Real) \quad y_\xi-1, U, U_\xi, \br, \br_\xi, \bs, H_\xi\in L^2(\Real) \cap L^\infty(\Real).
\end{equation}
These global solutions in Lagrangian coordinates can then be mapped to global
weak conservative solutions of the M2CH system as in \cite{holdenraynaud2007} in
the case of the CH equation. To be more specific, one has for each fixed time
$t$, cf. \cite[Definition 4.4]{guantanwei2015}, that
 \begin{equation}\label{map}
 u(x)=U(\xi), \quad \brho(x)=\br(\xi), \quad  \text{ and }\quad \brho_x(x)=\bs(\xi)
 \end{equation}
 for any $\xi$, such that $x=y(\xi)$ and 
 \begin{equation*}
 \mu(B)=\int_{\{ x\in y^{-1}(B)\}} H_\xi(\xi)d\xi \quad \text{ for any  Borel set } B,
 \end{equation*}
 where $\mu$ denotes the energy distribution measure. Since $\bar s(t,\cdot)\in
 L^\infty(\Real)$, it follows from \eqref{map} that $\brho_x$ remains also
 bounded in $L^\infty(\Real)$. In particular, it means that the blow-up of the
 solution only occurs when $u_x$ becomes unbounded, as in the CH case, the
 additional variable $\brho$ of M2CH does not blow up. Before closing this
 section, we introduce the \textit{Lagrangian potential energy} $r$ as
\begin{equation}
  \label{eq:defr}
  r(t,\xi) = \rho(t,y(t,\xi))y_\xi(t,\xi).
\end{equation}
As opposed to all the other Lagrangian variables introduced until now $(y
-\xi, U, H, \br, \bs)$, the Lagrangian variable $r$ is not generally bounded in
$L^\infty(\Real)$. Formally, $\br$ can be obtained from $r$ as
\begin{equation}\label{eq:rexp}
  \br(t,\xi)=\frac12 \int_\Real e^{-\abs{y(t,\xi)-y(t,\eta)}} r(t,\eta)d\eta
\end{equation}
and we also have the following relation between $\br$, $\bs$ and $r$, 
\begin{equation}
  \label{eq:relbsbrr}
  r = -\bs_\xi + \br y_\xi,
\end{equation}
from \eqref{eq:relrbrho}. From the definition of $r$, \eqref{eq:defr}, and the
transport equation \eqref{eq:m2chfull2}, we expect 
\begin{equation}
  \label{eq:goveeqforr}
  r_t = 0   
\end{equation}
This result can also be derived directly from the equivalent system
\eqref{eq:laggoveq} in purely Lagrangian variables. Indeed, after
differentiating \eqref{eq:relbsbrr} with respect to time, we get
\begin{equation*}
  r_t = - \bs_{\xi t} + \br_t y_\xi + \br y_{\xi t}.
\end{equation*}
We use \eqref{eq:laggoveq} and obtain
\begin{equation*}
  r_t = +(V_\xi + S_\xi) - (R + W)y_\xi + \br U_\xi.
\end{equation*}
From \eqref{eq:derVW} and \eqref{eq:derRS}, we obtain as expected that $r_t =
0$.
 
\section{Relabeling symmetry and local invariants}
\label{sec:relabsym}

\subsection{The case of the scalar Camassa--Holm equation}

As we mentioned in the introduction, the CH equation can be derived
from a variational principle, see \cite{constantin2002geometric} for a more
thorough presentation. In the case of the CH equation, there is no
potential energy and the Lagrangian is given by the kinetic energy only,
\begin{equation*}
  \LcalCH = \frac12 \int_\Real (u^2 + u_x^2)(t,x)\,dx.
\end{equation*}
The variation has to be done with respect to the particle path. We follow the
notations from \cite{constantin2002geometric} and denote the particle path by
$\phi(t,\xi)$, instead of $y(t,\xi)$ as in the previous section. After a change
of variable, we can rewrite $\LcalCH$ as
\begin{equation}
  \label{eq:LcalCH}
  \LcalCH(\phi) = \frac12\int_\Real \left(\phi_t^2\phi_\xi + \frac{\phi_{t\xi}^2}{\phi_\xi}\right)(t,\xi)\,d\xi.
\end{equation}
The group of diffeomorphism on $\Real$ lets the Lagrangian invariant with
respect to the group action of \textit{relabeling}. For a given diffeomorphism
$f$, the relabeling transformation of $\phi(t,\xi)$, with respect to $f$, is
given by $\phi\circ f = \phi(t, f(\xi))$. We can check directly that
\begin{equation*}
  \LcalCH(\phi\circ f)  = \frac12 \int_\Real \left((\phi_t\circ f)^2(\phi_\xi\circ f)f_\xi + \frac{(\phi_{t\xi}\circ f) ^2 f_\xi^2}{(\phi_\xi\circ f) f_\xi}\right)(t,\xi)\,d\xi =  \LcalCH(\phi),
\end{equation*}
after a change of variable. Noether's theorem tells us that to every
one-dimensional symmetry group which let the Lagrangian invariant, there
corresponds a conservation law. For the group of diffeomorphism, the tangent
space is formally isomorphic to $C^\infty(\Real)$, which is of infinite
dimension, so that we expect infinitely many invariants. Let us first shortly
present the Noether's argument in a finite dimensional setting, that is how a
symmetry leads to an invariant. We consider $q\in\Real^n$ and the Lagrangian
$\Lcal(q,\dot q)$. We assume that $\Lcal$ admits a one-dimensional symmetry
group. Keeping this presentation informal, we simply assume that there exists a
smooth mapping $S:\Real\times\Real^n\to\Real^n$, which represents the
one-dimensional group action, such that $S(0,\cdot) = \id$, and we denote
$q_\epsi(t)=S(\epsi, q(t))$. The invariance of the Lagrangian takes the form
\begin{equation}
  \label{eq:invfindim}
  \Lcal(q_\epsi(t), \dot q_\epsi(t)) = \Lcal(q(t), \dot q(t)).
\end{equation}
The Euler-Lagrange equations for the solution are 
\begin{equation}
  \label{eq:elfindim}
  \frac{d}{dt}\left(\fracpar{\Lcal}{\dot q}\right) = \fracpar{\Lcal}{q}.
\end{equation}
We differentiate \eqref{eq:invfindim} with respect to $\epsi$ and obtain
\begin{equation*}
  \fracpar{\Lcal}{q}\fracpar{q_\epsi}{\epsi} + \fracpar{\Lcal}{\dot q}\fracpar{\dot q_\epsi}{\epsi} = 0.
\end{equation*}
We set $\epsi=0$ and use the Euler-Lagrange equation and the previous equation
to get
\begin{equation*}
  \frac{d}{dt}\left(\fracpar{\Lcal}{\dot q}\right)\fracpar{q_\epsi}{\epsi}_{|\epsi = 0} + \fracpar{\Lcal}{\dot q}\fracpar{\dot q_\epsi}{\epsi}_{|\epsi = 0} = 0.
\end{equation*}
Assuming the solution is smooth, we have 
\begin{equation*}
  \fracpar{\dot q}{\epsi}_{|\epsi = 0} = \frac{d}{dt}\left(\fracpar{q}{\epsi}_{|\epsi = 0}\right)
\end{equation*}
and therefore it follows that
\begin{equation}
  \label{eq:pointconslawgen}
  \frac{d}{dt}\left(\fracpar{\Lcal}{\dot q}\fracpar{q_\epsi}{\epsi}_{|\epsi = 0}\right) = 0,
\end{equation}
and the quantity
\begin{equation*}
  \fracpar{\Lcal}{\dot q}\fracpar{q_\epsi}{\epsi}_{|\epsi = 0}
\end{equation*}
is preserved.  

Let us consider now the Lagrangian $\LcalCH$. To simplify the
notations we will denote the operator $\fracpar{}{\epsi}_{|\epsi = 0}$ as
$\delta_\epsi$. The Lagrangian is invariant with respect to
relabeling. Formally, the tangent space at the identity of the group of smooth
diffeomorphisms is the space of smooth functions $C^{\infty}(\Real)$. For any
function $g\in C^\infty(\Real)$ in the tangent space, we define the
diffeomorphism $f_\epsi(\xi) = \xi + \epsi g(\xi)$ and consider the
one-dimensional action defined as
\begin{equation*}
  \phi_\epsi = \phi\circ f_\epsi.
\end{equation*}
Slightly abusing the notation, we redefine $\LcalCH$ as
\begin{equation*}
  \LcalCH(\phi, \psi) = \frac12  \int_\Real \left(\psi^2\phi_\xi + \frac{\psi_{\xi}^2}{\phi_\xi}\right)(t,\xi)\,d\xi,
\end{equation*}
so that $\LcalCH(\phi, \partial_t\phi)$ is equal to $\LcalCH(\phi)$, as
introduced in \eqref{eq:LcalCH}.  The invariance of $\LcalCH$ with respect to
relabeling implies
\begin{equation*}
  \LcalCH(\phi_\epsi, \partial_t\phi_\epsi) =  \LcalCH(\phi, \partial_t\phi),
\end{equation*}
for all $\epsi\in\Real$, small enough. Following the same steps as before in the
finite dimensional case, we end up with the following conservation law, corresponding to
\eqref{eq:pointconslawgen},
\begin{equation}
  \label{eq:ptinvCH}
  \fracpar{}{t}\brac{\fracdelta{\LcalCH}{\psi},\delta_\epsi\phi}  = 0.
\end{equation}
Let us give a precise meaning to each of the expressions entering
\eqref{eq:ptinvCH}. We have
\begin{equation*}
  \brac{\fracdelta{\LcalCH}{\psi}, \delta\psi} = \int_\Real(\psi\delta\psi\phi_\xi + \frac{\psi_\xi}{\phi_\xi}\delta\psi_\xi) \,d\xi
\end{equation*}
so that 
\begin{equation*}
  \fracdelta{\LcalCH}{\psi}= \psi\phi_\xi  - \left(\frac{\psi_\xi}{\phi_\xi}\right)_\xi.
\end{equation*}
Let us compute $\delta_\epsi\phi_\epsi$. We have 
\begin{equation*}
  \delta_\epsi \phi_\epsi = \fracpar{}{\epsi}_{|\epsi = 0}\phi(\xi + \epsi g(\xi)) = \phi_\xi(\xi)g(\xi).
\end{equation*}
Hence, \eqref{eq:ptinvCH} can be rewritten as
\begin{equation*}
  \fracpar{}{t}\left(\int_\Real\left(\phi_{t}\phi_\xi  - \left(\frac{\phi_{t\xi}}{\phi_\xi}\right)_\xi\right)\phi_\xi g\,d\xi\right) = 0.
\end{equation*}
Assuming that the solution is smooth and decays sufficiently fast, it follows that
\begin{equation}
  \label{eq:variwithg}
  \int_\Real\fracpar{}{t}\left(\left(\phi_{t}\phi_\xi  - \left(\frac{\phi_{t\xi}}{\phi_\xi}\right)_\xi\right)\phi_\xi\right) g\,d\xi = 0,
\end{equation}
as $g$ is independent of time. Now, we use the fact that \eqref{eq:variwithg}
must hold for any $g\in C^\infty$ and therefore we obtain the following
\textit{pointwise} invariant,
\begin{equation*}
  \fracpar{}{t}\left(\left(\left(\phi_{t}\phi_\xi  - \left(\frac{\phi_{t\xi}}{\phi_\xi}\right)_\xi\right)\phi_\xi\right)(t,\xi)\right) = 0,
\end{equation*}
for all $\xi\in\Real$. Using the fact that $u\circ \phi = \phi_t$, $u_x\circ\phi
= \frac{\phi_{t\xi}}{\phi_\xi}$ and $u_{xx}\circ\phi
=\frac{(u_x\circ\phi)_\xi}{\phi_\xi}$, we can rewrite the pointwise invariant
above in the form of
\begin{equation}
  \label{eq:ptinvniceCH}
  \fracpar{}{t}((u - u_{xx})(t,\phi(t,\xi))\phi_{\xi}^2(t,\xi)) = 0.
\end{equation}
Note that, in fact, the pointwise invariant equation \eqref{eq:ptinvniceCH} can
be used to derive the CH equation in a rather straightforward way. To
see that, let us denote $m=\Helm u$ and expand \eqref{eq:ptinvniceCH}. We obtain
\begin{equation*}
  (m_t\circ\phi + m_x\circ\phi \phi_t)\phi_\xi^2 + 2m\circ\phi \phi_\xi\phi_{t\xi} = 0,
\end{equation*}
which, after using $\phi_t = u\circ\phi$ and $\phi_{t\xi} = u_x\circ\phi
\phi_\xi$, yields
\begin{equation*}
  (m_t + um_x + 2mu_x)\circ\phi \phi_\xi^2 = 0,
\end{equation*}
which, whenever $\phi_\xi$ does not vanish, is equivalent to the CH
equation. 

Multipeakon solutions, which are a special class of solutions for the CH
equation, are of the form
\begin{equation}
  \label{eq:mpCHdef}
  u(t,x)=\sum_{i=1}^n p_i(t)e^{-\vert x-q_i(t)\vert}
\end{equation}
for time-dependent coefficients $q_i(t)$, which denote the position of the
peaks, $p_i(t)$. Such solutions were identified in the seminal paper of Camassa
and Holm \cite{camassaholm1993}. Here, we want to show how the existence of this
special class of solutions can be inferred directly from the pointwise invariant
\eqref{eq:ptinvniceCH}. As the pointwise invariant is a consequence of the
symmetry of the problem, we can therefore establish a rather direct connection
between the symmetry of the Lagrangian and the existence of special
solutions. This approach is described in detail in \cite{holdenraynaud2007b} and
we sketch it here as a preparation for the case of M2CH. After applying the
Helmholtz operator $\Helm$ to $u$ in \eqref{eq:mpCHdef}, we get
\begin{equation}
  \label{eq:multpeakCHdef}
  (u - u_{xx})(t, x) = \sum_{i=1}^N2p_i(t)\delta(x - q_i(t)),
\end{equation}
where $\delta(x)$ denotes the Dirac delta distribution. For some initial data
that satisfies \eqref{eq:mpCHdef} and any point $x\in\Real$ away from the
singularities, meaning that it does not coincide with any of the $q_i$, we have
$(u - u_{xx} )(x)= 0$. After denoting $x(t)=y(t,\xi)$ the characteristic
starting at $x$, the pointwise invariant \eqref{eq:ptinvniceCH} yields $(u -
u_{xx})(t, x(t))=0$, as long as $y_\xi(t,\xi)\neq 0$. Hence, the structure given
by \eqref{eq:multpeakCHdef}, which defines the multipeakons, is preserved. The
formulation given by \eqref{eq:multpeakCHdef} cannot handle the collision of peaks
as some of the coefficients $p_i$ tend to $\pm\infty$ in this case. To handle
such case, we have to switch to the Lagrangian formulation. The pointwise
conservation equation plays then an essential role when showing that the multipeakon
structure is preserved. Between two neighboring peaks, say $q_i(t)$ and
$q_{i+1}(t)$, we have to show that $(u-u_{xx})(t, x) = 0$ for all $x\in(q_i(t),
q_{i+1}(t))$. The peaks follow the characteristics so that, in Lagrangian
coordinates, the region between the two peaks given as $\{(t,x)\ |\
q_i(t)<x<q_{i+1}(t)\}$, which is curved in Eulerian coordinates, becomes
\textit{rectangular}, that is $\{(t,\xi)\ |\ \xi_i < \xi < \xi_{i+1}\}$. Once
the pointwise conservation equation is established for each of such regions, we
can then deduce that the solution is indeed a multipeakon solution. The rigorous
presentation of this approach is given in \cite{holdenraynaud2007b}. 

\subsection{The case of the modified system}

The Lagrangian for the M2CH system is given in \eqref{eq:defLm2ch}. Let us
rewrite the potential energy in terms of the Lagrangian variables we have
introduced. We have 
\begin{equation*}
   \Epot = \frac12\int_\Real (\brho - \brho_{xx})\brho\,dx = \frac12\int_\Real \rho(x)\brho(x)\,dx.
\end{equation*}
We change to Lagrangian variable and obtain
\begin{equation*}
  \Epot = \frac12\int_\Real \br(\xi) r(\xi) \,d\xi.
\end{equation*}
We use the expression for $\br$ derived in \eqref{eq:rexp} and get
\begin{equation*}
  \Epot = \frac14\int_{\Real^2} e^{-\abs{\phi(t,\xi) - \phi(t,\eta)}} r(t,\eta)r(t,\xi)\,d\eta d\xi.
\end{equation*}
The relabeling transformation for a density variable such as $r$ is defined as
$r\mapsto r\circ f f_\xi$ for any $f\in\diff$. For such transformations, we can
check that the potential energy $\Epot$ is invariant. In \cite{grunert2016},
when we proceed with the variation for the 2CH system, the density
$\rho$ is treated as a function of $\phi$ so that the variation with respect to
$\rho$ is not computed independently. Here, we use a different approach by
decoupling the variables and introducing a Lagrangian multiplier function
$\lambda$ to enforce the mass conservation, that is $r_t=0$. Let $X = (\phi, r,
\lambda)$, we consider the Lagrangian defined as
\begin{multline}
  \label{eq:newLagm2ch}
  \Lcal(X, \partial_t X) = \frac12\int_\Real \left(\phi_t^2\phi_\xi +
    \frac{\phi_{t\xi}^2}{\phi_\xi}\right)(t,\xi)\,d\xi \\ - \frac14\int_{\Real^2}
  e^{-\abs{\phi(t,\xi) - \phi(t,\eta)}} r(t,\eta)r(t,\xi)\,d\eta d\xi\\ -
  \int_\Real \lambda(t,\xi) r_t(t,\xi)\,d\xi.
\end{multline}
We derive the Euler-Lagrange equation for this Lagrangian. Computations which we
only sketch here give us
\begin{subequations}
  \label{eq:dernewlag}
  \begin{multline}
    \label{eq:dernewlag1a}
    \fracdelta{\Lcal}{\phi} = -\frac12 (\phi_t^2)_\xi +\frac12
    \left(\frac{\phi_{t,\xi}^2}{\phi_\xi^2}\right)_\xi\\ +
    \frac{r}2\int_{\Real}\sign(\phi(t,\xi) - \phi(t,\eta)) e^{-\abs{\phi(t,\xi)
        - \phi(t,\eta)}} r(t,\eta)\,d\eta
  \end{multline}
  with
  \begin{equation}
    \label{eq:dernewlag1b}
    \fracdelta{\Lcal}{[\phi_t]} = \phi_t\phi_\xi  - \left(\frac{\phi_{t,\xi}}{\phi_\xi}\right)_\xi
  \end{equation}
  and
  \begin{align}
    \label{eq:dernewlag2}
    \fracdelta{\Lcal}{r} &= -\frac12\int_\Real
                                   e^{-\abs{\phi(t,\xi) - \phi(t,\eta)}} r(t,\eta)\,d\eta, &\fracdelta{\Lcal}{[r_t]}&=-\lambda,\\
    \label{eq:dernewlag3}
    \fracdelta{\Lcal}{\lambda} &= -r_t, &\fracdelta{\Lcal}{[\lambda_t]}&=0.
  \end{align}
\end{subequations}
We consider a diffeomorphism $\phi$ and a perturbation $\delta\phi$, then
$r(t,\xi)$ is perturbed by a corresponding $\delta r(t,\xi)$. Thus varying the
integral defined
\begin{equation}
\Lpot= \frac14 \int_{\Real^2} e^{-\vert \phi(t,\xi)-\phi(t,\eta)\vert } r(t,\eta) r(t,\xi)d\eta d\xi,
\end{equation}
with respect to $r$ yields
\begin{align*}
  \brac{\fracdelta{\Lpot}{r}, \delta r}& = \frac14 \int_{\Real^2} e^{-\vert \phi(t,\xi)-\phi(t,\eta)\vert } (\delta r(t,\eta) r(t,\xi)+r(t,\eta)\delta r(t,\xi)) d\eta d\xi\\
                                       & = \frac12 \int_{\Real}\left(\int_\Real e^{-\vert \phi(t,\xi)-\phi(t,\eta)\vert }r(t,\eta) d\eta \right)\delta r(t,\xi)d\xi,
\end{align*}
since we can interchange the order of integration.  Varying $\Lpot$ with respect
to $\phi$ yields
\begin{align*}
  &\Big<\fracdelta{\Lpot}{\phi},\delta\phi\Big>\\
       & = -\frac14 \int_{\Real^2} \sign(\phi(t,\xi)-\phi(t,\eta))e^{-\vert \phi(t,\xi)-\phi(t,\eta)\vert }r(t,\eta)r(t,\xi) (\delta\phi(t,\xi)-\delta\phi(t,\eta)) d\eta d\xi\\
       & = -\frac14 \int_{\Real}\int_\Real \sign(\phi(t,\xi)-\phi(t,\eta))e^{-\vert \phi(t,\xi)-\phi(t,\eta)\vert}r(t,\eta)d\eta r(t,\xi)\delta\phi(t,\xi)d\xi\\ 
       & \qquad +\frac14 \int_\Real \int_\Real \sign(\phi(t,\xi)-\phi(t,\eta)) e^{-\vert \phi(t,\xi)-\phi(t,\eta)\vert }  r(t,\eta) \delta\phi(t,\eta)d\eta r(t,\xi)d\xi\\
       & = -\frac12 \int_{\Real}\left(\int_\Real \sign(\phi(t,\xi)-\phi(t,\eta))e^{-\vert \phi(t,\xi)-\phi(t,\eta)\vert}r(t,\eta)d\eta r(t,\xi)\right)\delta\phi(t,\xi)d\xi,
\end{align*}
since we can again interchange the order of integration. The Euler-Lagrange
equation
\begin{equation*}
  \frac{d}{dt}\left(\fracdelta{\Lcal}{[X_t]}\right) =
  \fracdelta{\Lcal}{X}
\end{equation*}
yields 
\begin{equation}
  \label{eq:rtzero}
  r_t = 0,
\end{equation}
from \eqref{eq:dernewlag3} and
\begin{equation}
  \label{eq:lambdagov}
  \lambda_t = \frac12\int_\Real e^{-\abs{\phi(t,\xi) - \phi(t,\eta)}} r(t,\eta)\,d\eta,
\end{equation}
from \eqref{eq:dernewlag2}. Using the variable $\br$ defined as in
\eqref{eq:rexp}, we rewrite \eqref{eq:lambdagov} as
\begin{equation}
  \label{eq:lambdagov2}
  \lambda_t = \br.
\end{equation}
We can also rewrite \eqref{eq:dernewlag1a} as 
\begin{equation}
  \label{eq:simplcalphi}
  \fracdelta{\Lcal}{\phi} = -\frac12(\phi_t^2)_\xi +
  \frac12\left(\frac{\phi_{t,\xi}^2}{\phi_\xi^2}\right)_\xi - r\frac{\br_\xi}{\phi_\xi}.
\end{equation}
Then, \eqref{eq:simplcalphi} and \eqref{eq:dernewlag1b} yields
\begin{equation}
  \label{eq:govphi}
  \frac{d}{dt}\left(\phi_t\phi_\xi  - \left(\frac{\phi_{t,\xi}}{\phi_\xi}\right)_\xi\right) =  -\frac12(\phi_t^2)_\xi +
  \frac12\left(\frac{\phi_{t,\xi}^2}{\phi_\xi^2}\right)_\xi - r\frac{\br_\xi}{\phi_\xi}.
\end{equation}
The variable $r$ has been introduced as a primary variable, but since $r_t = 0$,
its dynamic is trivial. Setting $\rho\circ\phi = \frac{r}{\phi_\xi}$, $r_t = 0$
implies that
\begin{equation*}
  (\rho_t + (u\rho)_x)\circ\phi\phi_\xi = 0,
\end{equation*}
so that $\rho$ is indeed the density, if it is initially set as such. Moreover
we have
\begin{equation*}
  \br = \frac12\int_\Real e^{-\abs{\phi(t,\xi) - \phi(t,\eta)}} r(t,\eta)\,d\eta = \brho\circ\phi.
\end{equation*}
Hence,
\begin{equation*}
  r\frac{\br_\xi}{\phi_\xi} = \rho\circ\phi \brho_x\circ \phi \phi_\xi.
\end{equation*}
After some computation, we can then see that \eqref{eq:govphi} is equivalent to
\begin{equation*}
  (m_t + um_x + 2mu_x + \rho\brho_x)\circ\phi\phi_\xi = 0,
\end{equation*}
that is \eqref{eq:m2chfull1}, when $\phi_\xi$ does not vanish. Let us now
consider the action of relabeling on $\Lcal$ and derive pointwise
invariants. The action of the group on the Lagrangian multiplier $\lambda$ is
given by $(\lambda, f)\mapsto\lambda\circ f$, for any diffeomorphism $f$. As in
the scalar case of the CH equation, we consider for any $g\in C^\infty(\Real)$
the one-dimensional subgroup $f_\epsi(\xi)=\xi + \epsi g(\xi)$ of
diffeomorphisms. Using the notations introduced previously, we get
\begin{equation}
  \label{eq:varextvars}
  \delta_\epsi\phi = \phi_\xi g,\quad\delta_\epsi r = r_\xi g + rg_\xi = (rg)_\xi,\quad   \delta_\epsi\lambda = \lambda_\xi g.
\end{equation}
The pointwise conservation law  \eqref{eq:ptinvCH} becomes
\begin{equation}
  \label{eq:ptinvCH2}
  \fracpar{}{t}\left(\fracdelta{\LcalCH}{[\phi_t]}\delta_\epsi\phi + \fracdelta{\LcalCH}{[r_t]}\delta_\epsi r + \fracdelta{\LcalCH}{[\lambda_t]}\delta_\epsi\lambda\right)  = 0
\end{equation}
in this case. Hence, using \eqref{eq:dernewlag1b}, \eqref{eq:dernewlag2},
\eqref{eq:dernewlag3} and \eqref{eq:varextvars}, we get
\begin{equation}
  \label{eq:consgen1}
  \fracpar{}{t}\left(\int_\Real \left(\phi_t\phi_\xi  - \left(\frac{\phi_{t,\xi}}{\phi_\xi}\right)_\xi\right)\phi_{\xi} g\,d\xi - \int_\Real\lambda(rg)_\xi\,d\xi\right) = 0.
\end{equation}
Assuming that the solution is smooth and decays sufficiently fast, we move the
time derivative under the integral. The first integral is the same as in the
scalar case. For the second one, we get, after integration by parts,
\begin{equation*}
  \fracpar{}{t}\left(\int_\Real\lambda(rg)_\xi\,d\xi\right) =  -\int_\Real(\lambda_\xi r)_{t} g \,d\xi  
\end{equation*}
using the fact that $r_t = 0$. Hence, \eqref{eq:consgen1} yields
\begin{equation*}
  \int_\Real (m\circ\phi\phi_\xi^2 + \lambda_{\xi}r)_t g = 0,
\end{equation*}
which must hold for any function $g$ so that the pointwise conservation law for
the M2CH system is given by
\begin{equation}
  \label{eq:ptconslawm2ch0}
  ((m\circ\phi + \frac{\lambda_{\xi}}{\phi_\xi}\rho\circ\phi)\phi_{\xi}^2)_t = 0,
\end{equation}
and the pointwise conserved quantity is 
\begin{equation}
  \label{eq:pointlocainvm2ch}
  (m\circ\phi + \frac{\lambda_{\xi}}{\phi_\xi}\rho\circ\phi)\phi_{\xi}^2.
\end{equation}
Again, as for the case of CH, we observe that M2CH can be derived from
\eqref{eq:pointlocainvm2ch} in a rather straightforward manner. Using that $r_t
= 0$ and the expression \eqref{eq:lambdagov2} for $\lambda_t$, we get
\begin{equation*}
  (\lambda_{\xi}r)_t = \br_\xi r = (\rho\brho_x)\circ\phi\,\phi_\xi^2,
\end{equation*}
so that \eqref{eq:ptconslawm2ch0} can be rewritten as
\begin{equation}
  \label{eq:ptconslawm2ch}
  (m\circ\phi\,\phi_\xi^2)_t =- (\rho\brho_x)\circ\phi\,\phi_\xi^2,
\end{equation}
and, as before for the scalar case, we can check that \eqref{eq:ptconslawm2ch}
implies \eqref{eq:m2chfull1}.

From the pointwise conservation law \eqref{eq:ptconslawm2ch0}, we can deduce the
existence of multipeakon solutions. These are solutions $(u(t,x),\brho(t,x))$ of
the form
\begin{equation}
  \label{eq:mpdef}
  u(t,x)=\sum_{i=1}^n p_i(t)e^{-\vert x-q_i(t)\vert} \quad \text{ and }\quad \brho(t,x)=\sum_{i=1}^n s_i(t)e^{-\vert x-q_i(t)\vert} 
\end{equation}
for time-dependent coefficients $q_i(t)$, which denote the position of the
peaks, $p_i(t)$ and $s_i(t)$. After applying the Helmholtz operator $\Helm$ to $u$
and $\brho$ in \eqref{eq:mpdef} , we get
\begin{equation}
  \label{eq:multpeakdef}
  (u - u_{xx})(t, x) = \sum_{i=1}^N2p_i(t)\delta(x - q_i(t)),\quad \rho(t, x) = \sum_{i=1}^N2s_i(t)\delta(x - q_i(t)),
\end{equation}
where $\delta(x)$ denotes the Dirac delta distribution. Let us consider some
initial data that satisfies \eqref{eq:multpeakdef} initially. For any point
$x\in\Real$ away from the singularities, that is different from any of the
$q_i$, we have $(u - u_{xx} )(x)= 0$ and $\rho(x)=0$. Let us denote
$x(t)=y(t,\xi)$ the characteristic starting at $x$. Since $r_t=0$, we get that
$\rho(t, x(t))y_\xi(t,\xi)=0$, that is $\rho(t, x(t))=0$, as long as
$y_\xi(t,\xi)\neq 0$. From the pointwise invariant \eqref{eq:ptconslawm2ch0}, we
infer that $(u - u_{xx})(t, x(t))=0$, as long as $y_\xi(t,\xi)\neq 0$. Hence,
the structure given by \eqref{eq:multpeakdef}, which defines the multipeakons,
is preserved. The formulation given by \eqref{eq:multpeakdef} cannot handle the
collision of peaks as some of the coefficients $p_i$ tend to $\pm\infty$ in this
case. To do so, we have to switch to the Lagrangian formulation. To show that
the multipeakon structure is preserved in the Lagrangian formulation, the
pointwise conservation equation plays again an essential role but, clearly as
the following computations show, the derivation is significantly less
tractable. Between two neighboring peaks, say $q_i(t)$ and $q_{i+1}(t)$, we have
to show that $(u-u_{xx})(t, x) = 0$ and $\rho(t,x) = 0$ for all $x\in(q_i(t),
q_{i+1}(t))$. The peaks follow the characteristics so that, in Lagrangian
coordinates, the region between the two peaks given as $\{(t,x)\ |\
q_i(t)<x<q_{i+1}(t)\}$, which is curved in Eulerian coordinates, becomes
rectangular $\{(t,\xi)\ |\ \xi_i < \xi < \xi_{i+1}\}$. Once the pointwise
conservation equation is established for each of such regions, we can then
deduce that the solution is indeed a multipeakon solution. The rigorous
presentation of this approach is given in \cite{holdenraynaud2007b} and we only
sketch here how we prove the local conservation equation in the Lagrangian
setting. For each rectangular region of the form defined above, we can prove
that higher regularity for the Lagrangian variables is preserved by the
governing equations, see \cite{holdenraynaud2007b}. Then, we can define the
following quantities
\begin{equation*}
  (u - u_{xx})(t,y) y_\xi^2 = Uy_\xi^2 - U_{\xi\xi} + \frac{y_{\xi\xi}}{y_\xi}U_\xi
\end{equation*}
and 
\begin{equation}
  \label{eq:defr2}
  r = -\bs_\xi + \br y_\xi.
\end{equation}
Note that both quantities require higher regularity of the variables (existence
of $U_{\xi\xi}$, $y_{\xi\xi}$, $\bs_\xi$). For symplicity, we assume that
$y_\xi$ is different from zero. This assumption can then be removed as in
\cite{holdenraynaud2007b}. The pointwise conservation equation will be
established in the Lagrangian setting if we can show that the quantity $M$
defined below remains equal to zero,
\begin{equation}
  \label{eq:pwconseqelag}
  \underbrace{\frac{d}{dt}(Uy_\xi^2 - U_{\xi\xi} + \frac{y_{\xi\xi}}{y_\xi}U_\xi) + \bar r_\xi r + \lambda r_t}_{M} = 0.
\end{equation}
We have seen at the end of Section \ref{sec:lagvar} that $r_t=0$ can be derived
directly from the governing equations \eqref{eq:laggoveq} in Lagrangian
variables. Combining \eqref{eq:Qxi} and \eqref{eq:enidlag}, we get
\begin{equation}
  \label{eq:Qrew}
  Q_\xi = -U^2y_\xi - \frac12\frac{U_\xi^2}{y_\xi} +Py_\xi- \frac12\br^2y_\xi  + \frac12\bs^2y_\xi.
\end{equation}
Now, using the governing equations \eqref{eq:laggoveq}, we get
\begin{equation*}
  M = \underbrace{(U y_\xi^2)_t}_{A} + \underbrace{Q_{\xi\xi}}_{B} + \frac{U_{\xi\xi}U_\xi}{y_{\xi}} - \underbrace{\frac{y_{\xi\xi}}{y_\xi}Q_\xi}_{C} - \frac{y_{\xi\xi}U_\xi^2}{y_\xi^2} + \br_\xi r
\end{equation*}
with
\begin{align*}
  A &= -Qy_\xi^2 + 2UU_\xi y_\xi,\\
  B &= -2UU_\xi y_\xi - U^2y_{\xi\xi} - \frac{U_{\xi\xi}U_\xi}{y_\xi} + \frac12\frac{U_\xi^2y_{\xi\xi}}{y_\xi^2} + Qy_\xi^2 + Py_{\xi\xi} \\
    &\quad -\br\br_\xi y_\xi -\frac{\br^2}{2}y_{\xi\xi} + \bs\bs_\xi y_\xi + \frac{\bs^2}{2}y_{\xi\xi},\\
  C &= -U^2y_{\xi\xi} - \frac12 \frac{U_\xi^2y_{\xi\xi}}{y_\xi^ 2} +Py_{\xi\xi} -\frac{\br^2}{2}y_{\xi\xi} +\frac{\bs^2}{2}y_{\xi\xi}.
\end{align*}
Hence,
\begin{equation*}
  M = \br_\xi r - \br\br_\xi y_\xi + \bs\bs_\xi y_\xi.
\end{equation*}
and \eqref{eq:pwconseqelag} follows from \eqref{eq:defr2}.

\section{double multipeakons}
\label{sec:doublemultipeak}

For the CH equation the so-called multipeakon solutions serve on
the one hand as an illustrating example of how solutions may behave and on the
other hand they are dense in the set of weak conservative
solutions \cite{HoldenRaynaud2006b,HoldenRaynaud2008p}. Since the M2CH system reduces in the case $\rho\equiv 0$ to the CH
equation, the aim of this section is to derive the time evolution of solutions
until wave breaking in the case of both $u(t,x)$ and $\bar\rho(t,x)$ being
multipeakons, which we will call from now on double multipeakons. 
That is we are searching for solutions of the form \eqref{eq:mpdef}, where the positions of the peaks, $q_i(t)$, satisfy
\begin{equation*}
-\infty< q_1(t)< \dots <q_n(t)<\infty.
\end{equation*}
In particular, both $u(t,\cdot)$ and $\bar\rho(t,\cdot)$ are not differentiable at the points $x=q_i(t)$ ($i=1,2,\dots,n$), and hence $(u(t,x), \brho(t,x))$ are going to satisfy the M2CH system in the weak sense.

As a first step we have to define what it means to be a local weak solution of
the M2CH system. Direct computations as in \cite{holdenraynaud2007b} for the CH
equation yield that \eqref{eq:m2chfull} can be rewritten as follows
\begin{align*}
u_t-u_{txx}& +3uu_x-2u_xu_{xx}-uu_{xxx}+\brho\brho_x-\brho_x\brho_{xx}\\
& = u_t-u_{txx}+\frac32 (u^2)_x+\frac12 (u_x^2)_x-\frac12 (u^2)_{xxx}+\frac12 (\brho^2)_x-\frac12 (\brho^2_x)_x=0.
\end{align*}
and 
\begin{align*}
\rho_t+(u\rho)_x& = \brho_t-\brho_{txx} +u\brho_x+u_x\brho-u_x\brho_{xx}-u\brho_{xxx}\\
& = \brho_t-\brho_{txx}+(u\brho)_x-(u\brho)_{xxx}+(u_x\brho)_{xx}+(u_x\brho_x)_x=0.
\end{align*}
Hence we have the following definition.
\begin{definition}
We say that $(u, \bar\rho)\in L^1_{loc}([0,T], H^1_{loc})\times L^1_{loc}([0,T], H^1_{loc})$ is a weak solution of the M2CH system if it satisfies 
\begin{subequations}\label{eq:weakm2CH}
  \begin{align}\label{eq:weakm2CH1}
    u_t-u_{txx}+\frac32 (u^2)_x+\frac12 (u_x^2)_x-\frac12 (u^2)_{xxx}+\frac12 (\brho^2)_x-\frac12 (\brho^2_x)_x&=0\\ \label{eq:weakm2CH2}
    \brho_t-\brho_{txx}+(u\brho)_x-(u\brho)_{xxx}+(u_x\brho)_{xx}+(u_x\brho_x)_x&=0
  \end{align}
\end{subequations}
in the sense of distributions.
\end{definition}
Since the local, weak multipeakon solutions are piecewise smooth solutions and
following closely the computations carried out in \cite{HoldenRaynaud2006} for
the CH equation, we obtain after some integration by parts that all the
information concerning the time evolution of $q_i(t)$, $p_i(t)$, and $s_i(t)$,
is contained in the coefficients of $\delta_{q_i}$ and $\delta'_{q_i}$. For
\eqref{eq:weakm2CH1} the coefficient of $\delta'_{q_i}$ must be equal to zero
and is given by
\begin{equation*}
  [u_t]_{q_i}+\frac12 [(u^2)_x]_{q_i}=0,
\end{equation*}
where we denote by $[v]_{q_i}=v_i(q_i+)-v_{i-1}(q_i-)$. The jumps are mainly influenced by the sign changes in the derivative, which come from the term $p_i(t)e^{-\vert x-q_i(t)\vert }$ at the point $x=q_i(t)$. In particular, we have 
\begin{equation*}
  [u_t]_{q_i}=2p_i(t)q_i'(t), \quad [(u^2)_x]_{q_i}=-4p_i(t)u(q_i(t)),
\end{equation*}
and hence
\begin{equation*}
  2p_i(t)(q_i'(t)-u(q_i(t)))=0.
\end{equation*}
Dividing both sides by $2p_i(t)$ yields the equation for the characteristic
\begin{equation*}
  q_i'(t)=u(q_i(t))=\sum_{j=1}^n p_j(t)e^{-\vert q_j(t)-q_i(t)\vert}.
\end{equation*}

By the same argument the coefficient of $\delta_{q_i}$ in \eqref{eq:weakm2CH1} must be equal to zero and thus
\begin{equation*}
  [u_{t,x}]_{q_i}-\frac12 [u_x^2]_{q_i}+\frac12 [(u^2)_{xx}]_{q_i}+\frac12 [\brho_x^2]_{q_i}=0.
\end{equation*}
Again the jumps are mainly influenced by the sign changes in the derivatives, which come from the terms $p_i(t)e^{-\vert x-q_i(t)\vert}$ and $s_i(t)e^{-\vert x-q_i(t)\vert}$, it is therefore convenient to introduce the following abbreviations
\begin{equation*}
  u(t,x)=\sum_{j=1}^n p_j(t)e^{-\vert x-q_j(t)\vert } = f(t,x)+p_i(t)e^{-\vert x-q_i(t)\vert }
\end{equation*}
and 
\begin{equation*}
  \brho(t,x)=\sum_{j=1}^n s_j(t)e^{-\vert x-q_j(t)\vert }=g(t,x)+s_i(t)e^{-\vert x-q_i(t)\vert}.
\end{equation*}
Direct computations, similar to the ones before, then yield
\begin{align*}
  [u_{t,x}]_{q_i}=-2p_i'(t), \quad & [u_x^2]_{q_i}=-4p_i(t)f_x(t, q_i(t)),\\
  [(u^2)_{xx}]_{q_i}= -8p_i(t)f_x(t, q_i(t)), \quad & [\brho_x^2]_{q_i}=-4s_i(t)g_x(t,q_i(t)),
\end{align*}
which implies that 
\begin{equation*}
  -2p_i'(t)+2p_i(t)f_x(t,q_i(t))-4p_i(t)f_x(t,q_i(t))-2s_i(t)g_x(t,q_i(t))=0.
\end{equation*}
Recalling the definition of $f(t,x)$ and $g(t,x)$, we end up with 
\begin{align*}
  p_i'(t)& =\sum_{j\not =i} p_i(t)p_j(t)\sign(q_i(t)-q_j(t))e^{-\vert q_i(t)-q_j(t)\vert}\\
         & \quad  +\sum_{j\not = i} s_i(t)s_j(t) \sign(q_i(t)-q_j(t))e^{-\vert q_i(t)-q_j(t)\vert }.
\end{align*}

As far as \eqref{eq:weakm2CH2} is concerned, the coefficient of $\delta'_{q_i}$ has to be equal to zero and is given by  
\begin{equation*}
  [\brho_t]_{q_i}+[(u\brho)_x]_{q_i}-[u_x\brho]_{q_i}=0.
\end{equation*}
In particular, we have 
\begin{align*}
  & [\brho_t]_{q_i}=2s_i(t)q_i'(t), \quad [u_x\brho]_{q_i}=-2p_i(t)(s_i(t)+g(t,q_i(t))),\\ 
  & [(u\brho)_x]_{q_i}= -2(s_i(t)f(t,q_i(t))+p_i(t)g(t,q_i(t))+2s_i(t)p_i(t)),
\end{align*} 
and accordingly
\begin{equation*}
  2s_i(t)q_i'(t)-2s_i(t)(f(t,q_i(t))+p_i(t))=0
\end{equation*}
Recalling the definition of $f(t,x)$ and dividing both sides by $s_i(t)$ we obtain
\begin{equation*}
  q_i'(t)=f(t,q_i(t))+p_i(t)=u(q_i(t)).
\end{equation*}

By the same argument the coefficient of $\delta_{q_i}$ must be equal to zero, which is equivalent to 
\begin{equation*}
  [\brho_{t,x}]_{q_i}+[(u\brho)_{xx}]_{q_i}-[(u_x\brho)_x]_{q_i}-[u_x\brho_x]_{q_i}=0.
\end{equation*}
Direct computations yield
\begin{align*}
  [\brho_{t,x}]_{q_i}& =-2s_i'(t),\\
  [(u\brho)_{xx}]_{q_i}& =-4p_i(t)g_x(t,q_i(t))-4s_i(t)f_x(t, q_i(t))\\
  [u_x\brho_x]_{q_i}& = -2p_i(t)g_x(t,q_i(t))-2s_i(t)f_x(t, q_i(t)),\\
  [(u_x\brho)_x]_{q_i}& =-2p_i(t)g_x(t,q_i(t))-2s_i(t)f_x(t, q_i(t)),
\end{align*}
and hence
\begin{align*}
  s_i'(t)=0.
\end{align*}

Thus we have the following system of ODEs
\begin{subequations}\label{sys:ODE}
  \begin{align}
    q_i'(t)&=\sum_{i=1}^n p_i(t)e^{-\vert q_i(t)-q_j(t)\vert },\\
    p_i'(t)& = \sum_{j\not =i} (p_i(t)p_j(t)+s_i(t)s_j(t))\sign(q_i(t)-q_j(t))e^{-\vert q_i(t)-q_j(t)\vert },\\
    s_i'(t)&=0.
  \end{align}
\end{subequations}

\section{Double peakon-antipeakon solutions}

In this section, we study in detail the peakon-antipeakon solutions in the case
$n=2$. That means both $u(t,x)$ and $\bar\rho(t,x)$ are the sum of one peakon
and one antipeakon except when wave breaking occurs, in which case both are
constantly equal to zero and part of the energy is concentrated in one point,
which is represented by a $\delta$-distribution. To set the stage, let
\begin{subequations}\label{ansatz}
  \begin{align}
    u(t,x)&= p_1(t)e^{-\vert x-q_1(t)\vert} +p_2(t)e^{-\vert x-q_2(t)\vert}, \\
    \bar\rho(t,x)&= s_1(t)e^{-\vert x-q_1(t)\vert}+s_2(t)e^{-\vert x-q_2(t)\vert}.
  \end{align}
\end{subequations}
We assume that $q_1\leq q_2$ initially and, as the peaks travel along
characteristics, this property remains true for all time. Then the corresponding
time independent total energy, which we denote $E$ is given by
\begin{equation}\label{Hamiltonian}
  E = p_1^2(t)+p_2^2(t)+s_1^2(t)+s_2^2(t)+2(p_1p_2(t)+s_1s_2(t))e^{q_1(t)-q_2(t)}.
\end{equation}
Introducing the variables $q=q_1-q_2$, $Q=q_1+q_2$, $p=p_1-p_2$, $P=p_1+p_2$,
$s=s_1-s_2$, and $S=s_1+s_2$, \eqref{ansatz} and \eqref{Hamiltonian} rewrite as
\begin{subequations}
  \begin{align}
    u(t,x)&=\frac12 (p+P)(t)e^{-\vert x-\frac12(q+Q)(t)\vert}+\frac12(P-p)(t)e^{-\vert x-\frac12(Q-q)(t)\vert},\\
    \bar\rho(t,x)& = \frac12 (s+S)(t)e^{-\vert x -\frac12 (q+Q)(t)\vert}+\frac12(S-s)(t)e^{-\vert x-\frac12 (Q-q)(t)\vert},\\
    \label{eq:Hmult}
    E&=\frac12 (p^2+P^2+s^2+S^2)(t)+\frac12 (P^2-p^2+S^2-s^2)(t)e^{q(t)}.
  \end{align}
\end{subequations}
According to \eqref{sys:ODE}, the functions $q$, $Q$, $p$, $P$, $s$, and $S$ satisfy the following system of ordinary differential equations
\begin{subequations}
  \begin{align}
    q_t(t)&=p(t)(1-e^{q(t)}),\quad & Q_t(t)&= P(t)(1+e^{q(t)}),\\
    p_t(t)&= \frac 12 p(t)^2+\frac12 C,\quad & P_t(t)&=0,\\
    s_t(t)&= 0, \quad &S_t(t)&= 0,
  \end{align}
\end{subequations}
where $C=(P^2(t)+S^2(t)+s^2(t)-2E)$. We observe that if $Q(t)=P(t)=S(t)=0$ holds
for some $t$, then it holds for all $t$. This means, since $Q(t)=q_1(t)+q_2(t)$,
$P(t)=p_1(t)+p_2(t)$, and $S(t)=s_1(t)+s_2(t)$, that there exist two peakon
solutions $(u(t,x),\brho(t,x))$ such that
\begin{subequations}
  \begin{align}
    u(t,x)&= p_1(t)(e^{-\vert x-q_1(t)\vert } - e^{-\vert x+q_1(t)\vert})\\
    \brho(t,x)&= s_1(t)(e^{-\vert x-q_1(t)\vert} - e^{-\vert x+q_1(t)\vert}).
  \end{align}
\end{subequations}
Such solutions are called peakon-antipeakon solutions, since both $u(t,\cdot)$
and $\brho(t,\cdot)$ are antisymmetric for all $t\in\Real$. In the remaining of
this section, we compute these solutions explicitly. Wave breaking occurs when
two peakons occupy the same position, that is $q(t^*)= q_1(t^*)-q_2(t^*)=0$. In
this case, we have $u_x(t,x)\to \mp \infty$ as $t\to t^*\mp$, which implies
$p(t)\to\pm\infty$ as $t\to t^*\mp$. As mentioned earlier, $\bar\rho$ and its
derivative $\bar\rho_x$ remain bounded.  We now turn to the computation of
$p(t)$, $q(t)$, and $u(t,q_1(t))$, the value of $u(t,x)$ at the left peak,
depending on the value of $s^2/2$ compared with the total energy $E$. We observe
that the governing equations are invariant with respect to the transformation
\begin{equation*}
  t\mapsto\alpha t,\quad u\mapsto\alpha u,\quad \brho\mapsto \alpha\brho.
\end{equation*}
Therefore, we do not restrict ourselves by considering only a single value of
$E$. For simplicity, we choose $E=\frac12$ so that \eqref{eq:Hmult} yields
\begin{equation}
  \label{eq:enpresone}
  (p(t)^2 + s(t)^2)(1 - e^{q(t)}) = 1.
\end{equation} 
Moreover, the equation is also invariant by the transformation
$\brho\mapsto-\brho$ so that, without loss of generality, we assume $s\geq 0$.
Let us denote by $u_\dagger$ and $\brho_\dagger$ the values of $u$ and $\brho$
at the peaks, that is
\begin{equation*}
  u_\dagger(t) = u(t, q_1(t)) = -u(t, q_2(t)) = \frac12 p(t)(1 - e^{q(t)})
\end{equation*}
and
\begin{equation*}
  \brho_\dagger(t) = \brho(t,q_1(t)) = -\brho(t,q_2(t))  = \frac12 s(t)(1 - e^{q(t)}).
\end{equation*}
From these expressions, we can express $s$ and $p$ as function of $u_\dagger$ and
$\brho_\dagger$ and plug the results in \eqref{eq:enpresone}. We obtain
\begin{equation*}
  4u_\dagger^2 + 4\brho_\dagger^2 = (1 - e^q).
\end{equation*}
We use again the definition of $\brho_\dagger$ and get
\begin{equation}
  \label{eq:trajcircle}
  u_\dagger^2 + \left(\brho_\dagger - \frac{1}{4s}\right)^2 = \left(\frac1{4s}\right)^2.
\end{equation}
Since $s$ is constant, the trajectories of $(\brho_\dagger, u_\dagger)$ lie on
circles as depicted in Figure \ref{fig:trajplot}. Let us know consider the
following three cases, depending on the value of $s$, which cover all the
possible types of dynamics for the system, and compute explicitly the solution
for each case with initial data $q(0)=q_0$ and $p(0)=p_0$.

\subsection{Case $0\leq s <1$}
In this case, we have
\begin{equation*}
  p_t=\frac12 (p-\sqrt{-C})(p+\sqrt{-C})
\end{equation*}
We integrate this expression and obtain
\begin{subequations}\label{eq:Formel2}
  \begin{align}
    p(t) &= \sqrt{-C}\frac{1+Ae^{\sqrt{-C}t}}{1-Ae^{\sqrt{-C}t}},   \\
    q(t) &= -\ln\left(1+e^{-\sqrt{-C}t}\frac{(1-Ae^{\sqrt{-C}t})^2}{(1-A)^2}(e^{-q_0}-1)\right)
  \end{align}
  where $A=\frac{p_0-\sqrt{-C}}{p_0+\sqrt{-C}}$ and
  \begin{equation*}
    u_\dagger(t) = \frac12 \sqrt{-C}e^{-\sqrt{-C}t} \frac{(1-Ae^{\sqrt{-C}t})(1+Ae^{\sqrt{-C}t})(e^{-q_0}-1)}{(1-A)^2 +e^{-\sqrt{-C}t}(1-Ae^{\sqrt{-C}t})^2 (e^{-q_0}-1)}.
  \end{equation*}
\end{subequations}
By definition, we have $C=s^2-1$ and $(p_0^2+s^2)(1-e^{q_0})=1$. Hence,
\begin{equation*}
  p_0^2=-C+(p_0^2+s_0^2)e^{q_0}\geq -C
\end{equation*}
and $A>0$ so that wave breaking occurs at time
$t^*=\frac{1}{\sqrt{-C}}\ln(\frac1{A})$. We shift time so that the collision
takes place at $t=0$. To do so, we let $p_0$ tend to infinity and $q_0$ to zero,
while preserving \eqref{eq:enpresone}, that is $(p_0^2+s^2)(1 - e^{q_0}) =
1$. \def\pinf{p_\infty} Let us denote
\begin{equation*}
  \pinf = \sqrt{1 - s^2}=\sqrt{-C}.  
\end{equation*}
The solution \eqref{eq:Formel2} is equivalent, up to shift in time, to
\begin{align*}
  p(t) &= \pinf\frac{1 + e^{\pinf t}}{{1 - e^{\pinf t}}},\\
  q(t) &= -\ln\left(1 + \frac{\cosh(\pinf t) - 1}{2\pinf^2}\right).
\end{align*}
Moreover, we have
\begin{equation*}
  u_\dagger(t) = -\pinf\frac{\sinh(\pinf t)}{4\pinf^2+2(\cosh(\pinf t) - 1)}.
\end{equation*}
When $t$ tends to $\pm\infty$, we get the following limits
\begin{equation*}
  \lim_{t\pm\infty}u_\dagger(t) = \mp\frac{\pinf}{2}.
\end{equation*}
Let us write $u_\dagger(\pm\infty)=\lim_{t\pm\infty} u_\dagger (t)$ and use the
same notation for $p(\pm\infty)$ and $\brho_\dagger(\pm\infty)$.  Taking the
same limit in \eqref{eq:enpresone}, we get $p(\pm\infty)^2 + s^2 = 1$, which
implies
\begin{equation}
  \label{eq:limcirc}
  u_\dagger(\pm\infty)^2 + \brho_\dagger(\pm\infty)^2 = \frac14.
\end{equation}
This circle is plotted in Figure \ref{fig:trajplot} and represents the limiting
values for $(\brho_\dagger, u_\dagger)$.

\subsection{Case $s > 1$}
The solution reads for all $t\in\Real$
\begin{subequations}
  \begin{align*}
    p(t) &= \sqrt{C}\tan\left(\frac12\sqrt{C}t+D\right),\\
    q(t) &= -\ln\left(1+B\cos^2(\frac12\sqrt{C}t+D)\right)
  \end{align*}
  and we have
  \begin{equation*}
    u_\dagger(t) = \frac14\sqrt{C}\frac{B\sin( \sqrt{C}t+2D)}{1+B\cos^2(\frac12 \sqrt{C}t+D)}  \end{equation*}
\end{subequations}
where
\begin{equation}
  D=\arctan\left(\frac{p_0}{\sqrt{C}}\right)\quad \text{ and }\quad B=\frac{e^{-q_0}-1}{\cos(D)^2}.
\end{equation}
We shift time as before and set the collision time to zero. The solution is then
given by
\begin{subequations}
  \begin{align*}
    p(t) &=- \sqrt{C}\cot\left(\sqrt{C}t/2\right),\\
    q(t) &= -\ln\big(1+\frac1{C}\sin^2\left(\sqrt{C}t/2\right)\big)
  \end{align*} 
  and we have
  \begin{align*}
    u_\dagger(t) &= -\frac14\sqrt{C}\frac{\frac1{C}\sin(\sqrt{C}t)}{1+\frac{1}{C}\sin^2(\sqrt{C}t/2)},\\ 
  \end{align*}
\end{subequations}
Especially the last double peakon-antipeakon solution comes as a surprise, since
such peakon-antipeakon solutions do not exist for the CH equation.  In the case
of the CH equation, i.e. $s=S=0$ for all $t\in\Real$, the constant $C$ reduces
to
\begin{equation*}
  C=-2=-p_0^2(1-e^{q_0}),
\end{equation*} 
which is less than $0$ under the assumption that $q_0\not=0$ and
$u(0,x)\not\equiv 0$. 
For the M2CH system, on the other hand, 
\begin{equation*}
  C=-p^2_0(1-e^{q_0})+s^2e^{q_0}.
\end{equation*}
and choosing $s$ big in contrast to $p_0$, one ends up in the case $C\geq0$. Thus the last case is intrinsic for the M2CH system.

\subsection{Case $s = 1$}

Direct calculations in that case yield
\begin{subequations}\label{Formel:1}
  \begin{align*}
    p(t) &=\frac{2p_0}{2-tp_0},\\  
    q(t) &= -\ln\left(1+\left(\frac{2-tp_0}{2}\right)^2(e^{-q_0}-1)\right)
  \end{align*}
  and
  \begin{equation*}
    u_\dagger(t) = \frac{p_0\frac{2-tp_0}{4}(e^{-q_0}-1)}{1+\left(\frac{2-tp_0}{2}\right)^2 (e^{-q_0}-1)},
  \end{equation*}
\end{subequations}
As in the previous case, we set the collision time to zero and obtain
\begin{subequations}
  \begin{align*}
    p(t) &=-\frac{1}{t/2},\\  
    q(t) &= -\ln\left(1+\left(t/2\right)^2\right)
  \end{align*}
  and
  \begin{align*}
    u_\dagger(t) &= -\frac{t}{t^2 + 4}.
  \end{align*}
\end{subequations}
In this limiting case $u_\dagger(\infty)=0$.

\begin{figure}[h]
  \centering
  \def\insertPlot #1{\includegraphics[width=5cm]{lowelas-0#1}}
  \begin{tabular}[t]{cc}
    \insertPlot{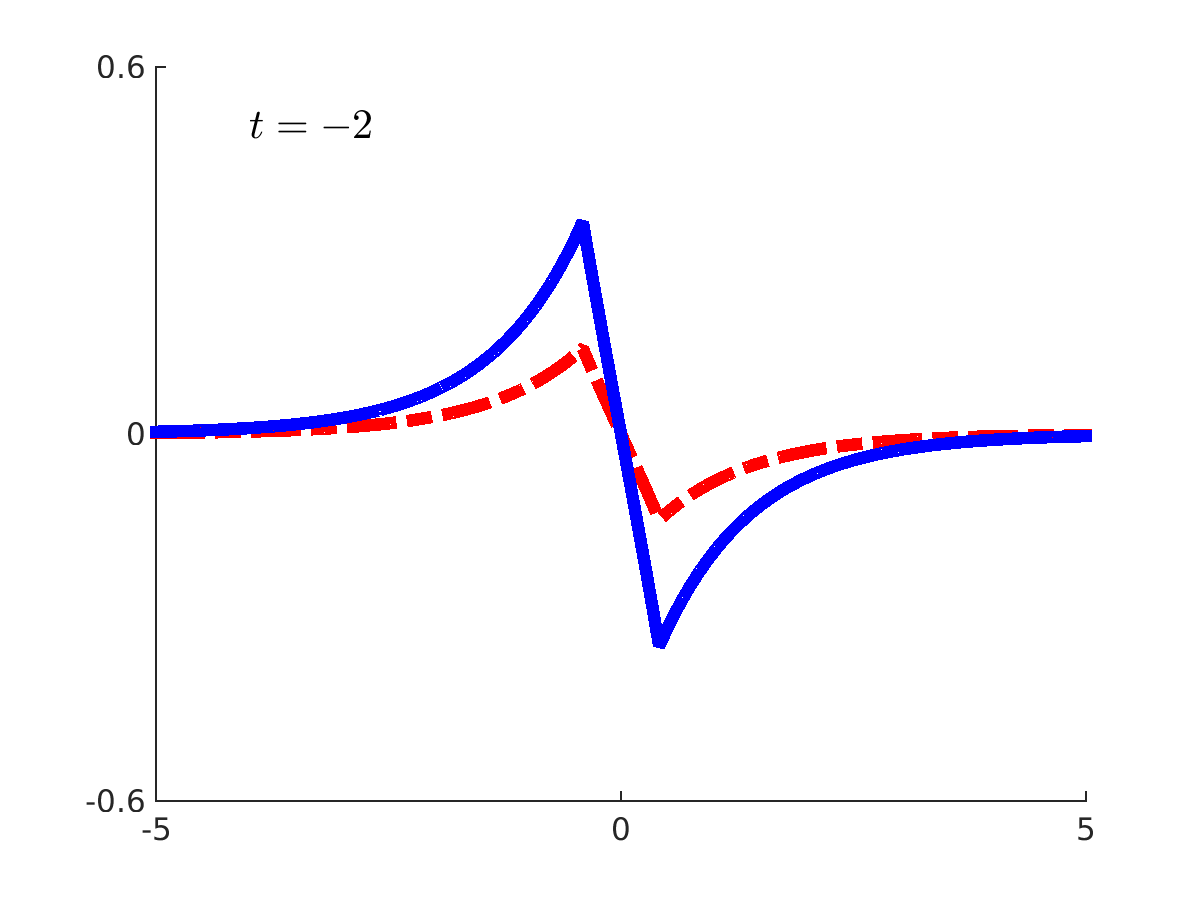}&\insertPlot{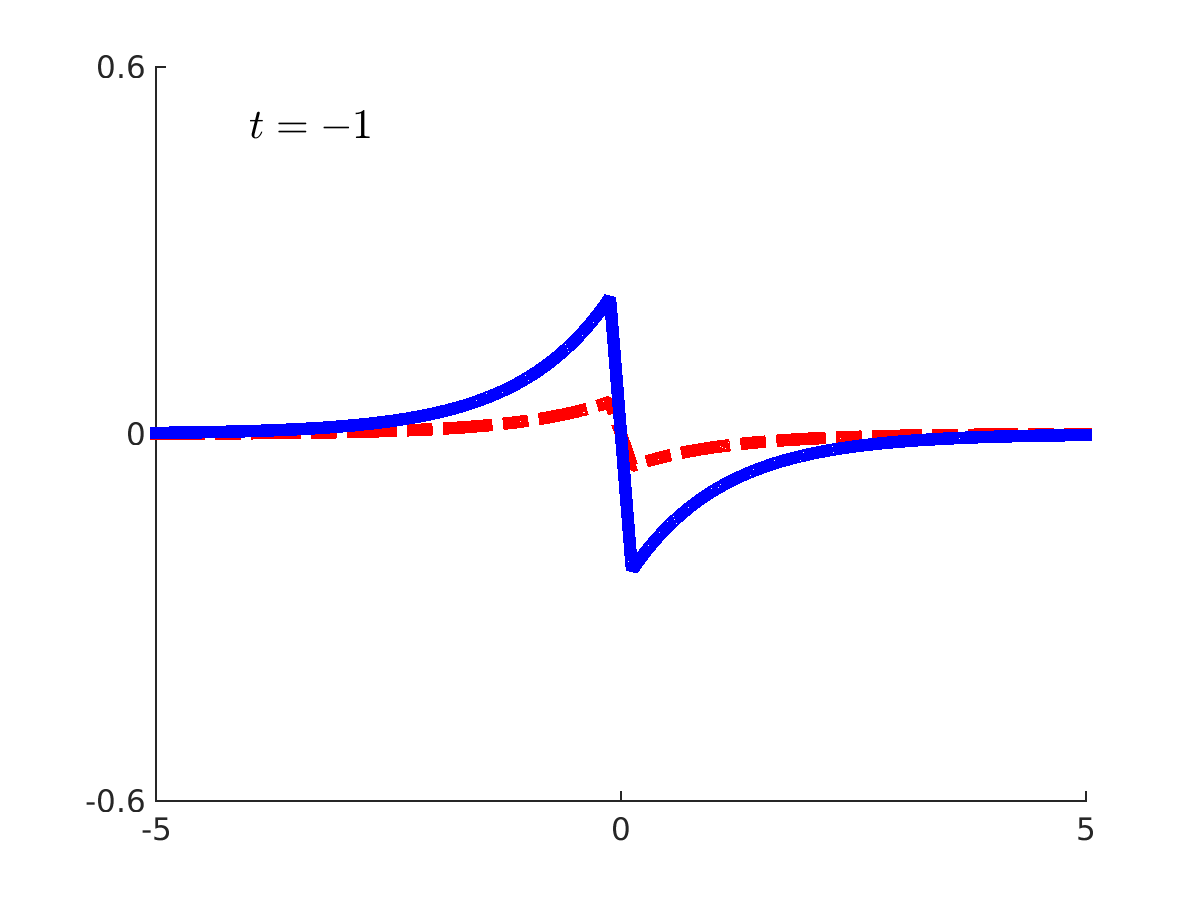}\\
    \insertPlot{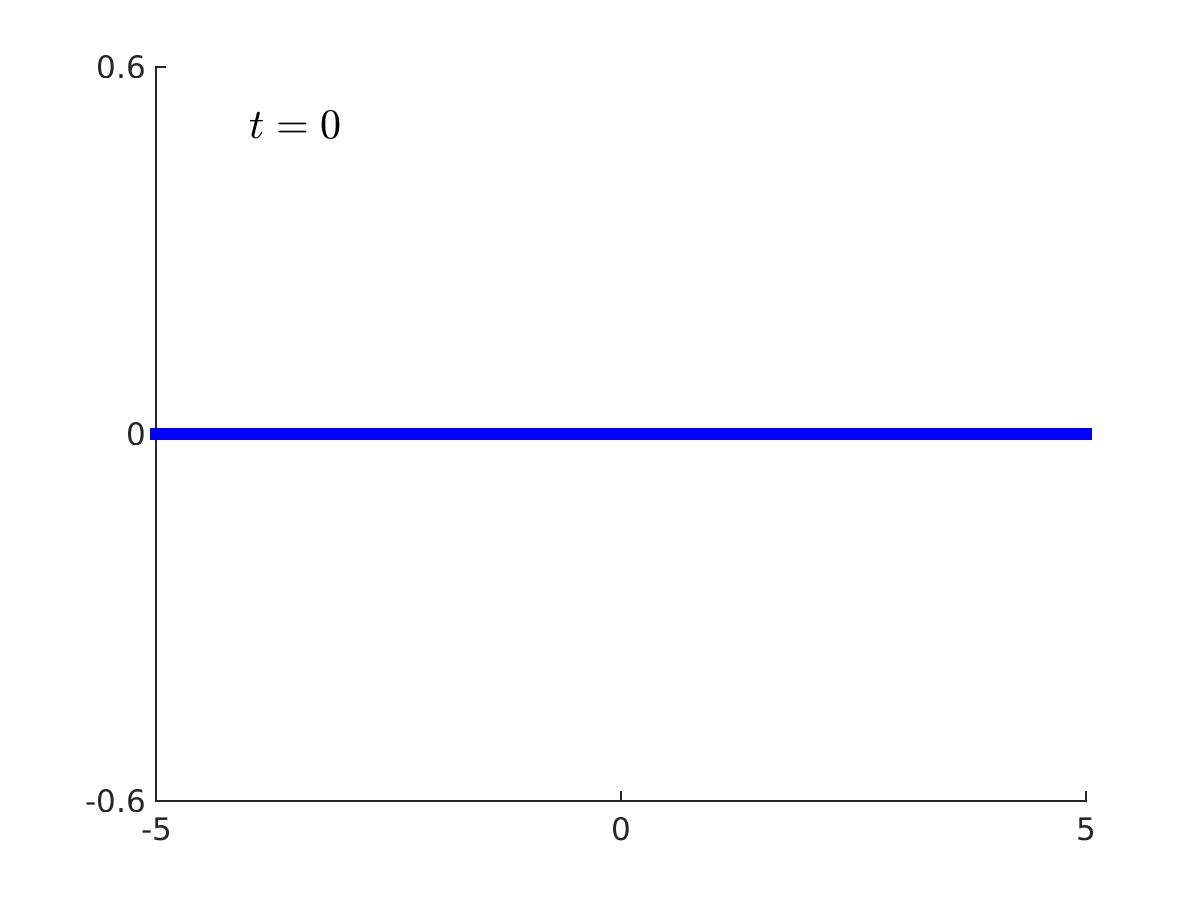}&\insertPlot{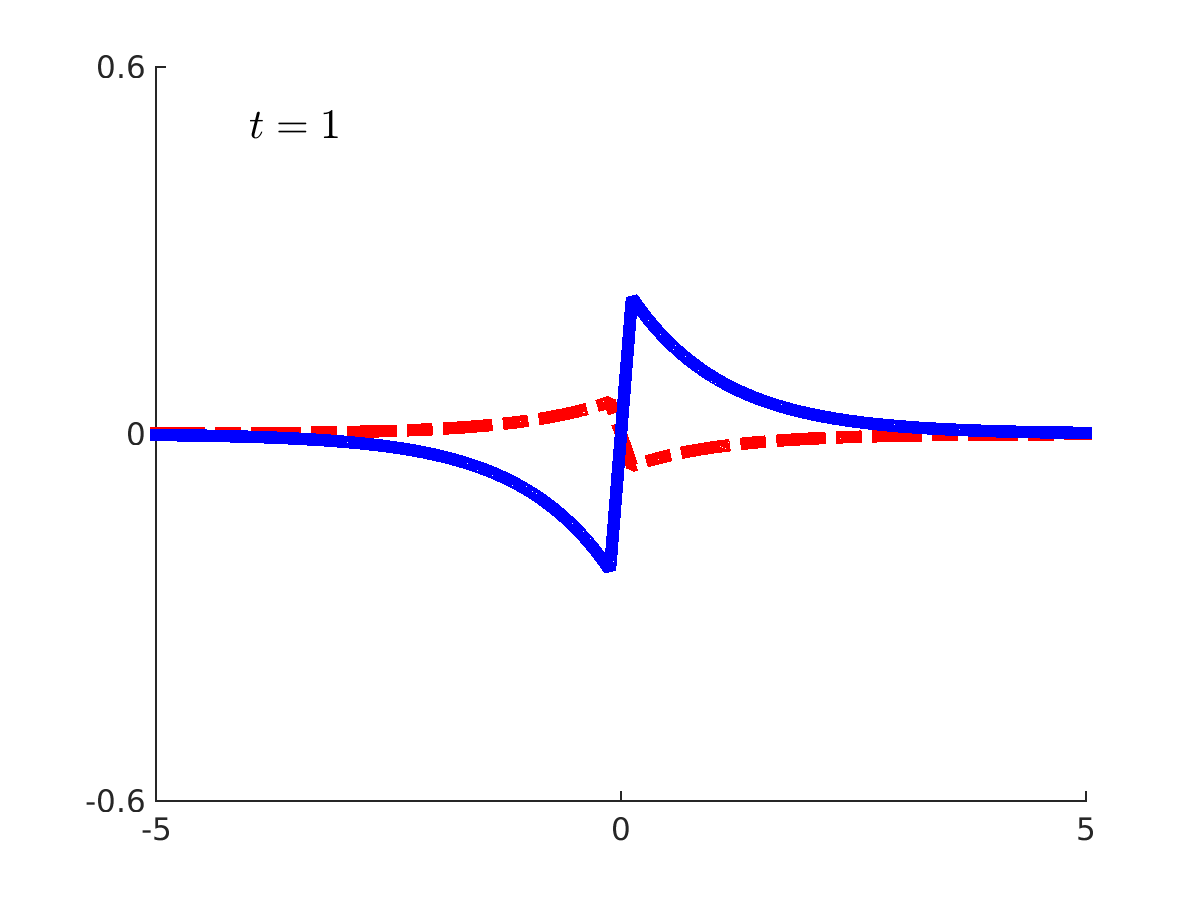}\\
    \insertPlot{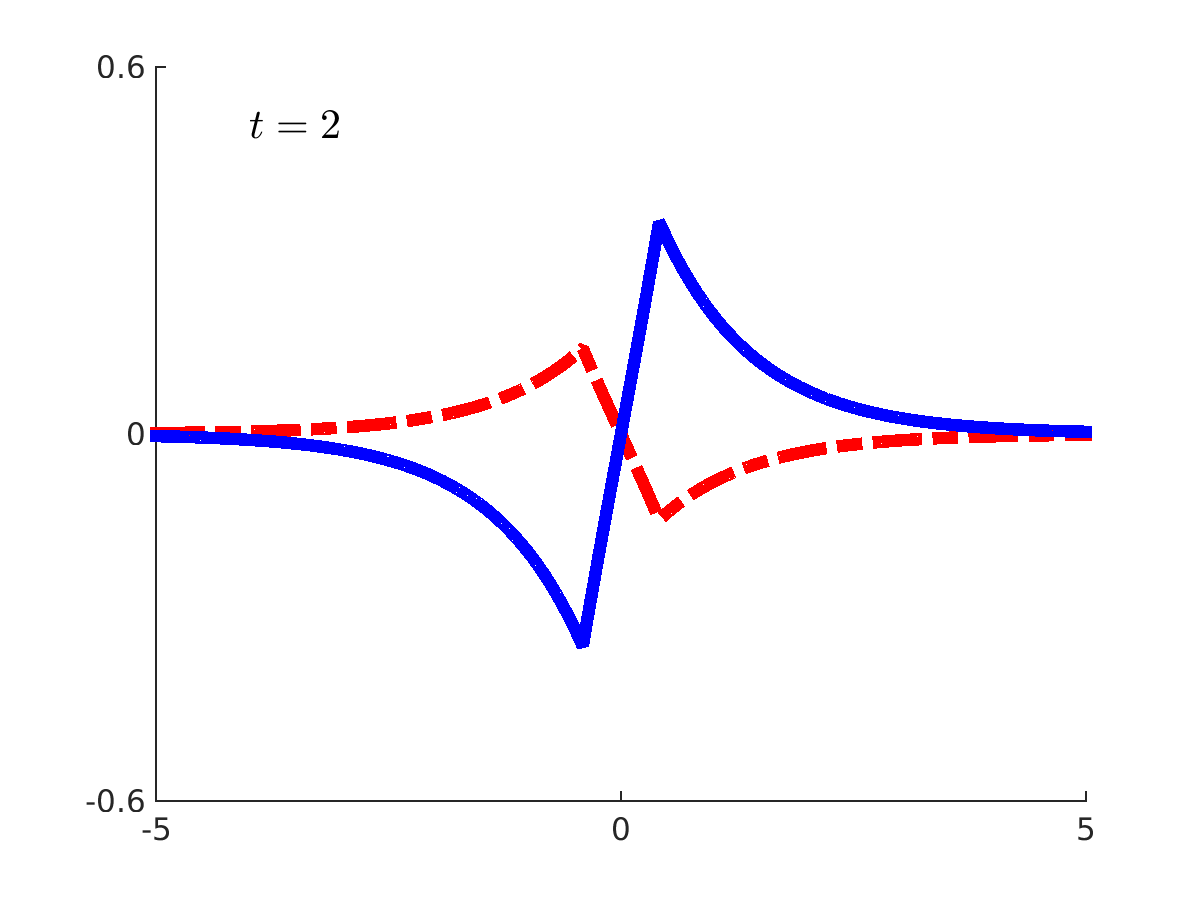}&\insertPlot{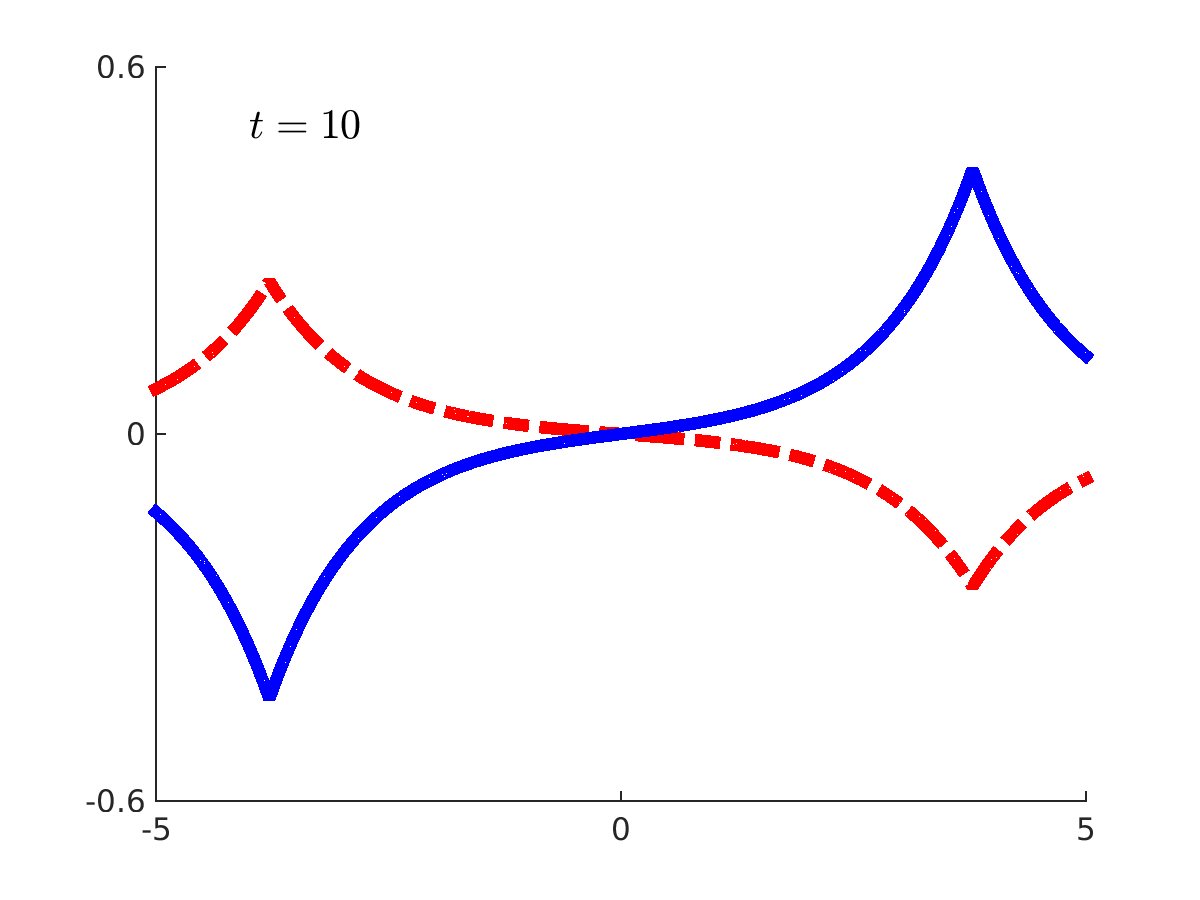}
  \end{tabular}
  \caption{Case $0 \leq s < 1$. Plot of the solution $u$ (blue) and $\brho$ (red) at
    different times.}
  \label{fig:lowelas}
\end{figure}

\begin{figure}[h]
  \centering
  \def\insertPlot #1{\includegraphics[width=5cm]{highelas-0#1}}
  \begin{tabular}[t]{cc}
    \insertPlot{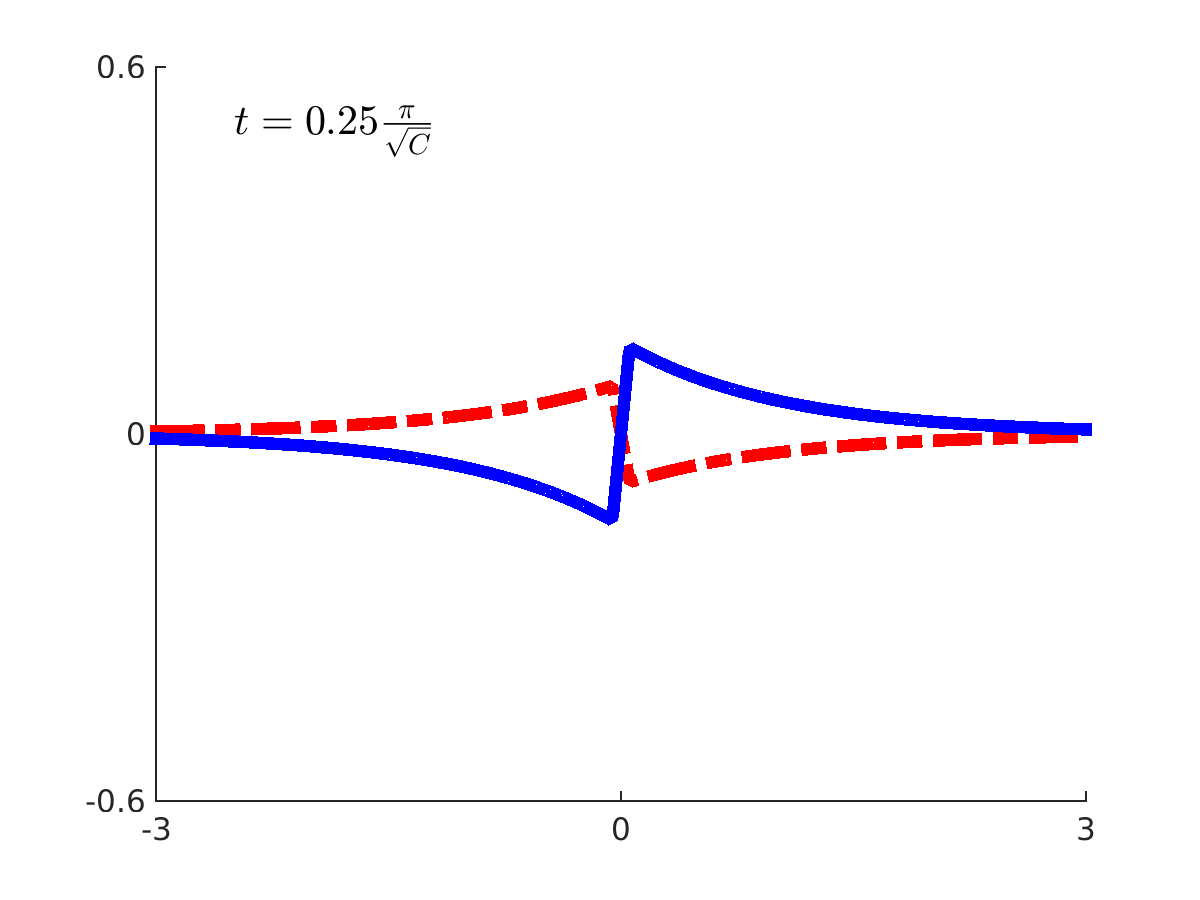}&\insertPlot{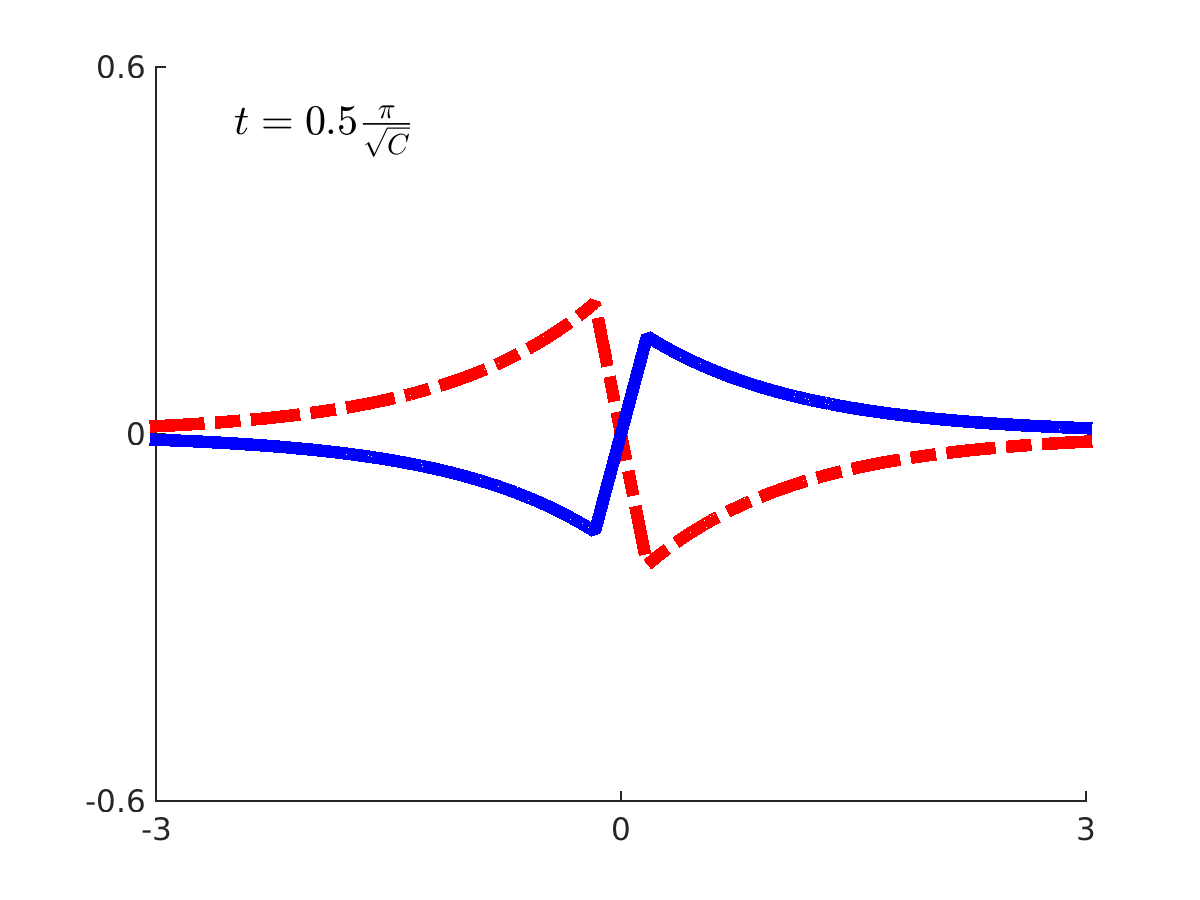}\\
    \insertPlot{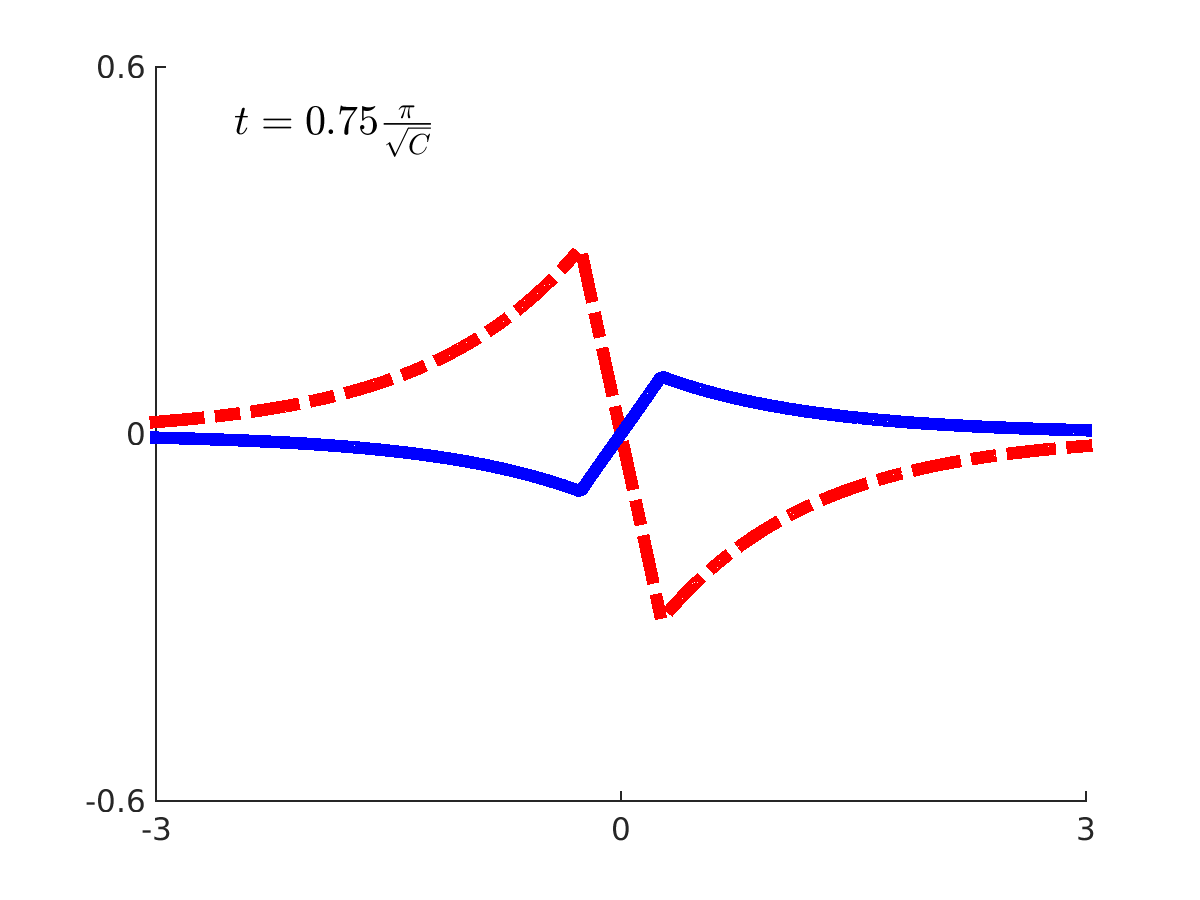}&\insertPlot{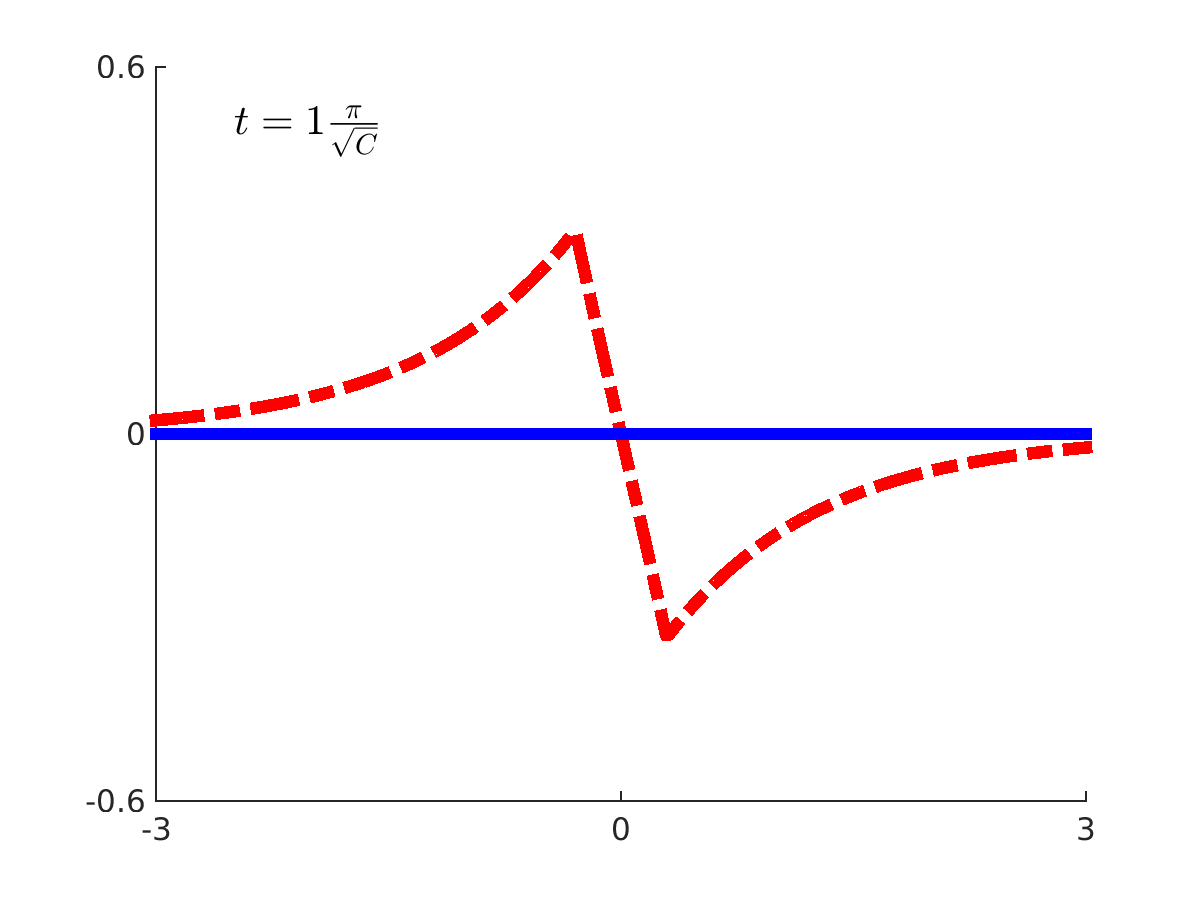}\\
    \insertPlot{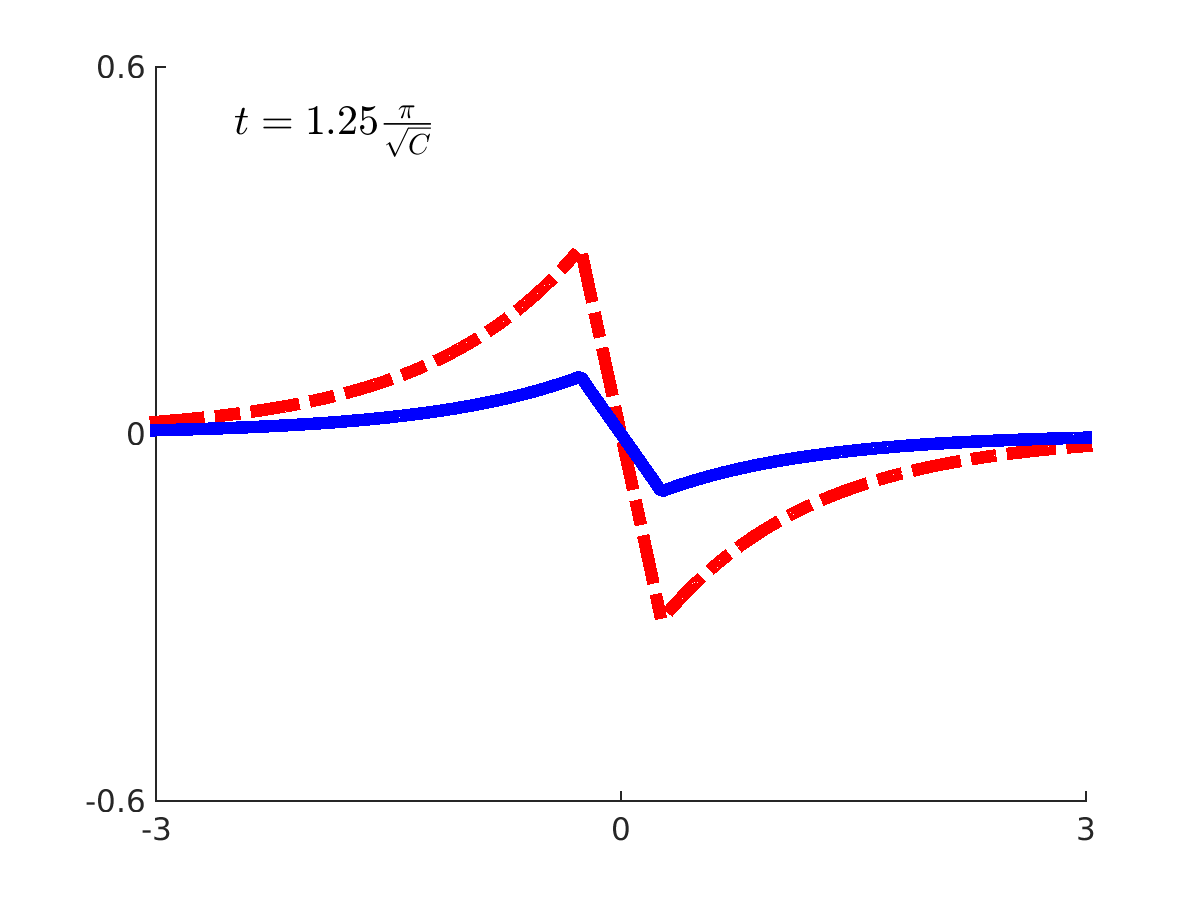}&\insertPlot{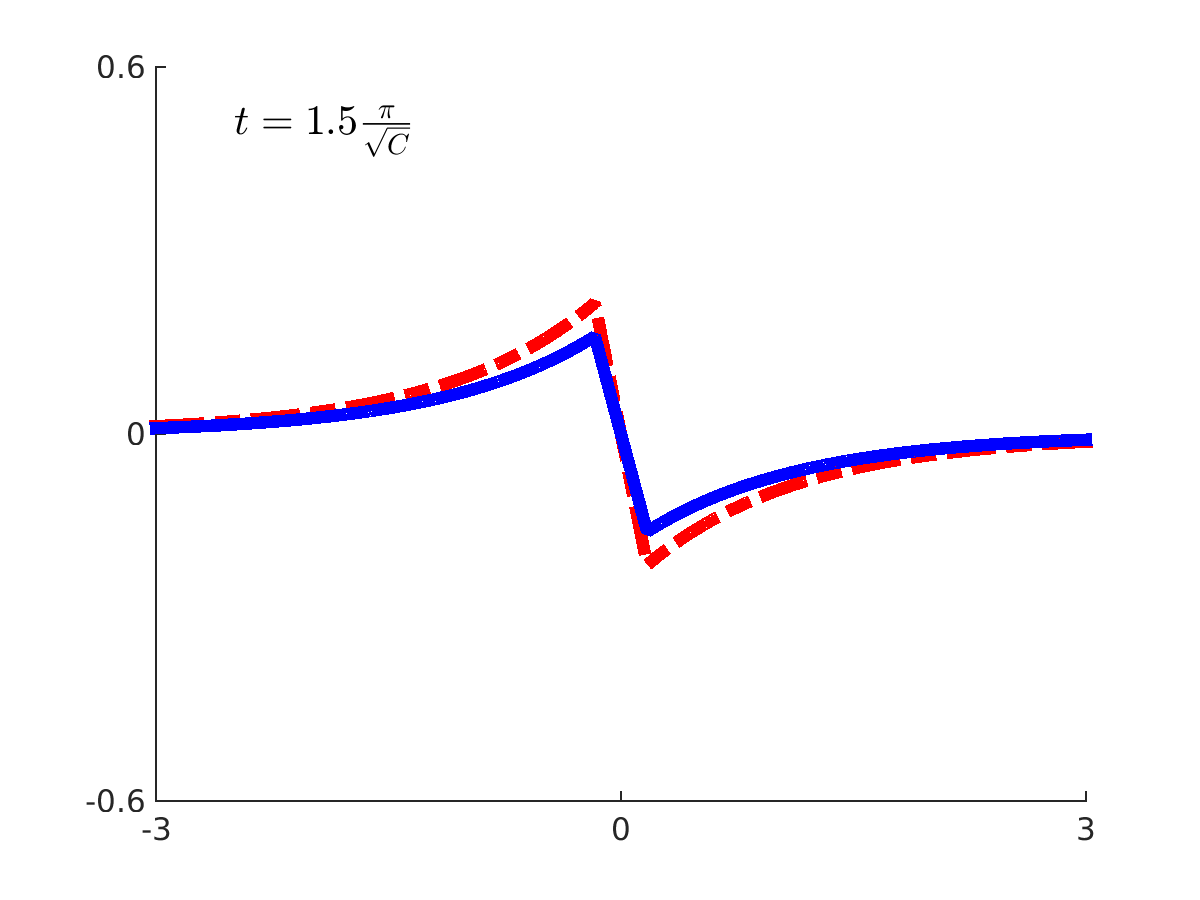}\\
    \insertPlot{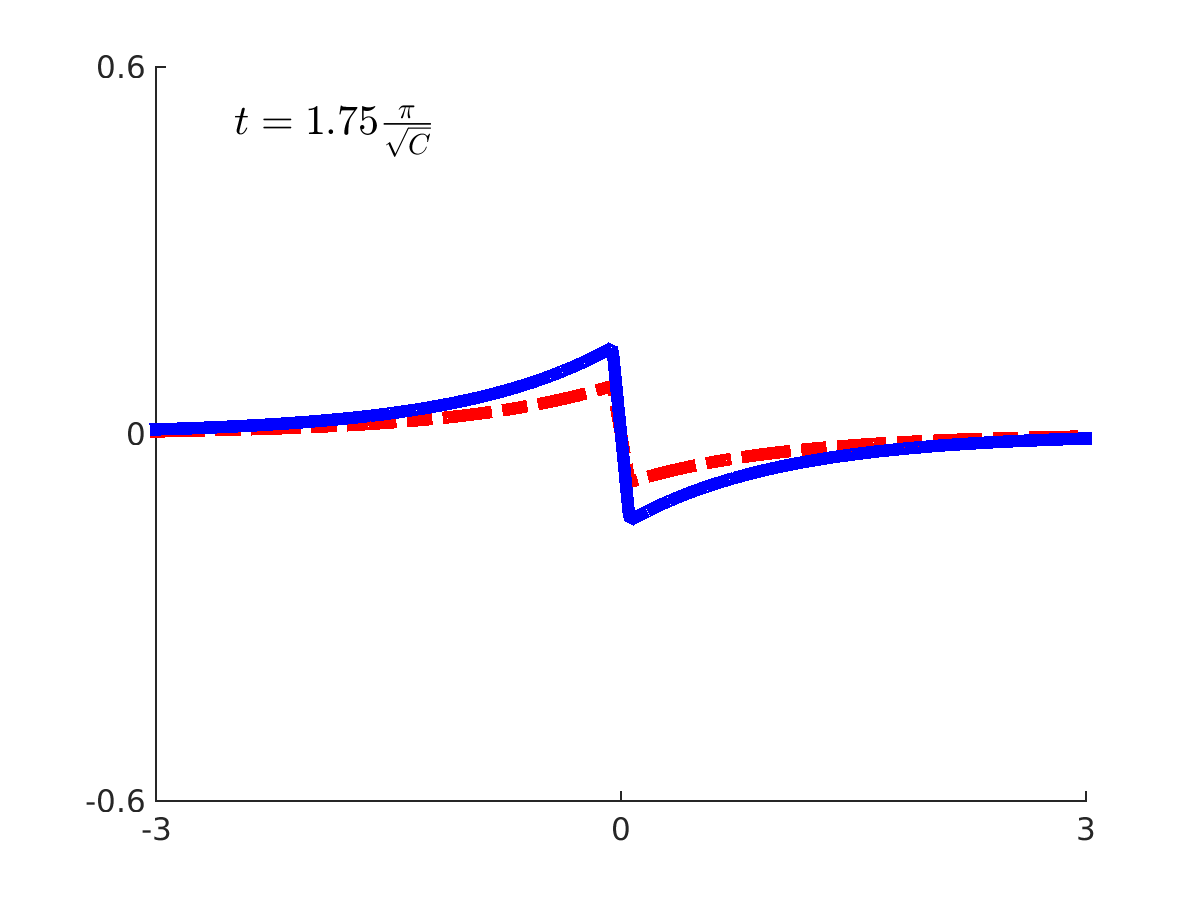}&\insertPlot{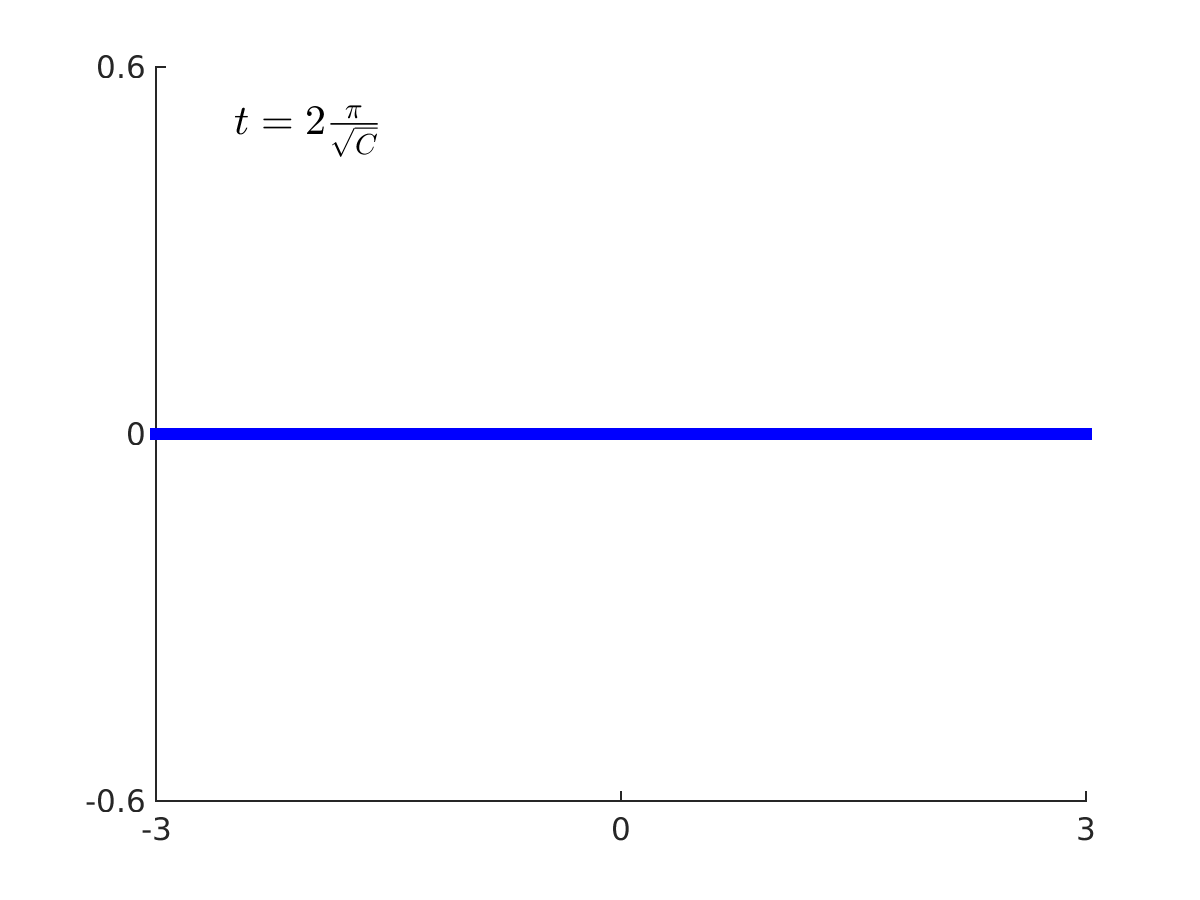}
  \end{tabular}
  \caption{Case $s > 1$. The solution is periodic with period
    $\frac{2\pi}{\sqrt{C}}$. The first plot (top, left) shows the solution right
    after a collision.}
  \label{fig:highelas}
\end{figure}

\begin{figure}[h]
  \centering
  \def\insertPlot #1{\includegraphics[width=5cm]{limcase-0#1}}
  \begin{tabular}[t]{cc}
    \insertPlot{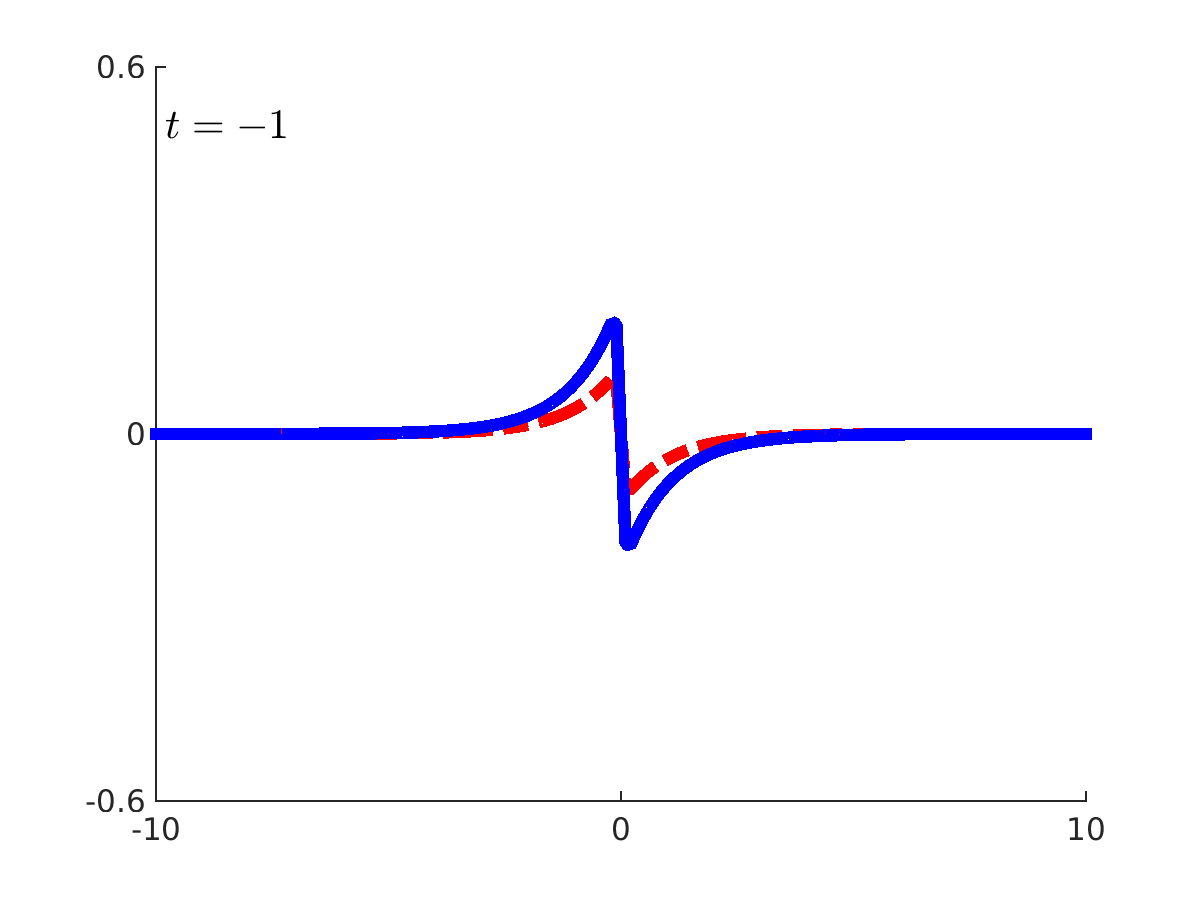}&\insertPlot{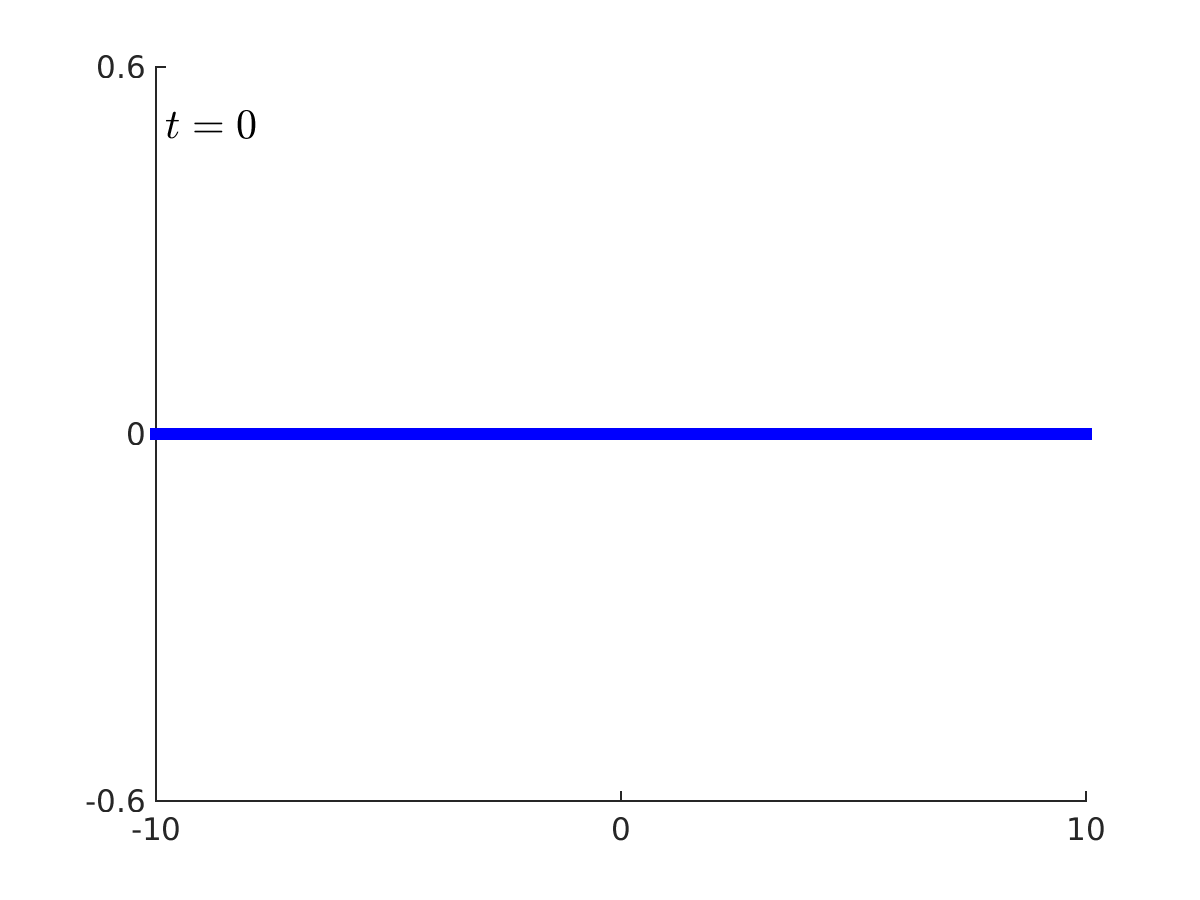}\\
    \insertPlot{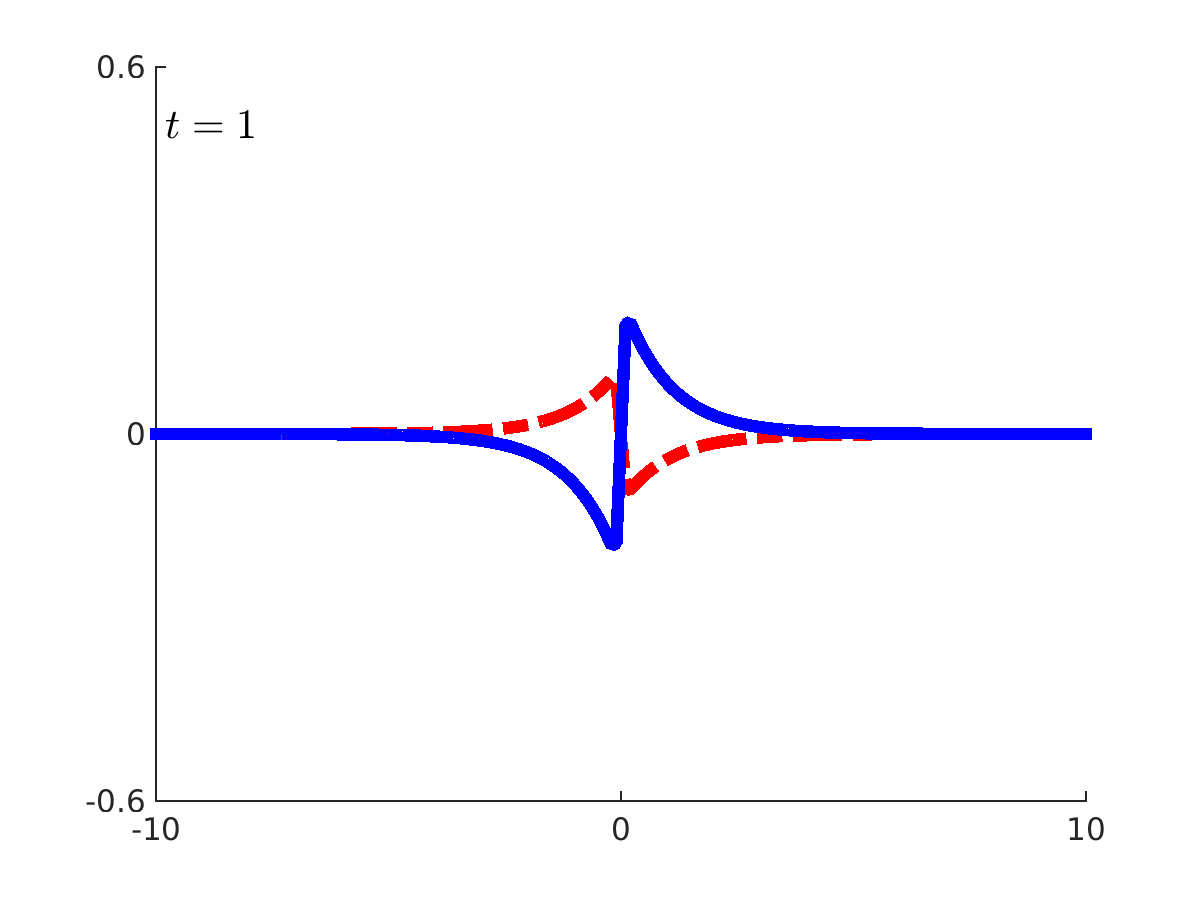}&\insertPlot{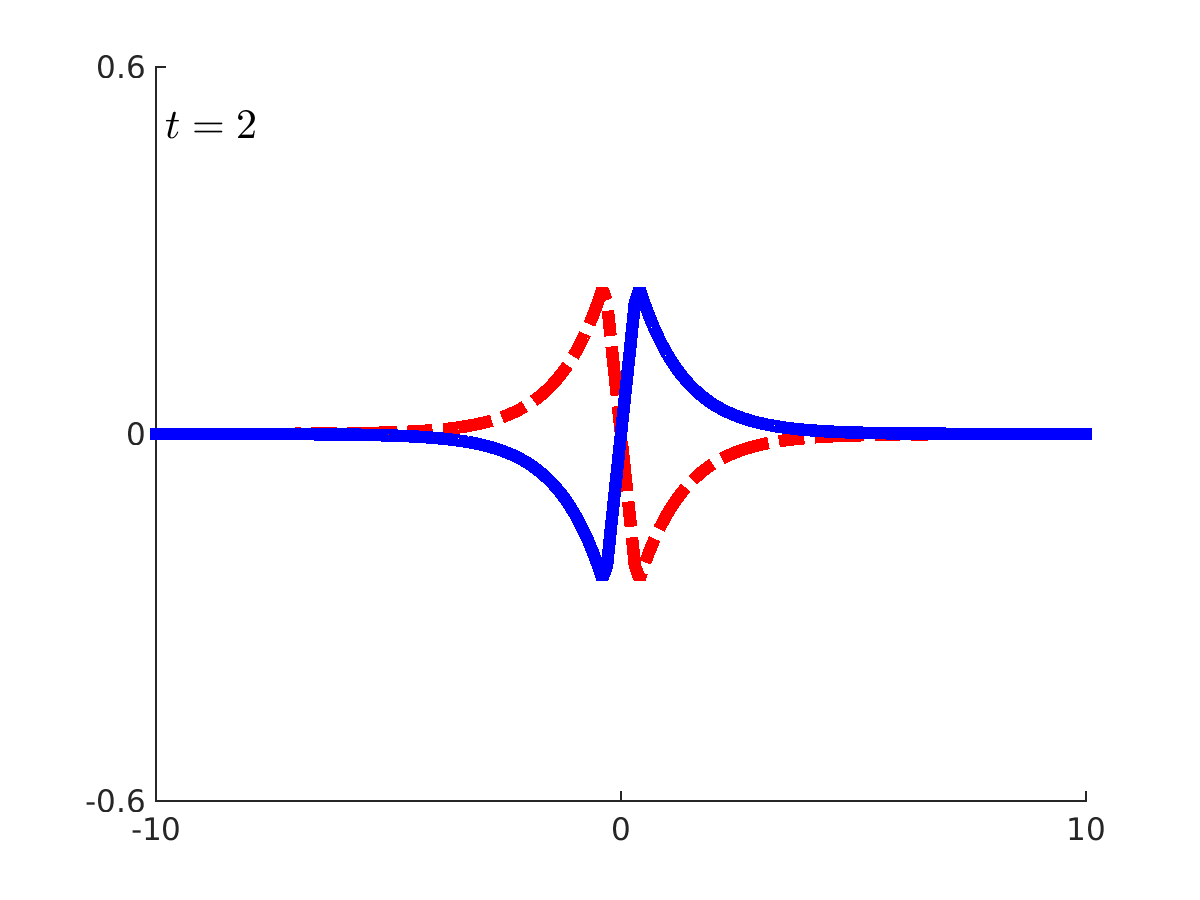}\\
    \insertPlot{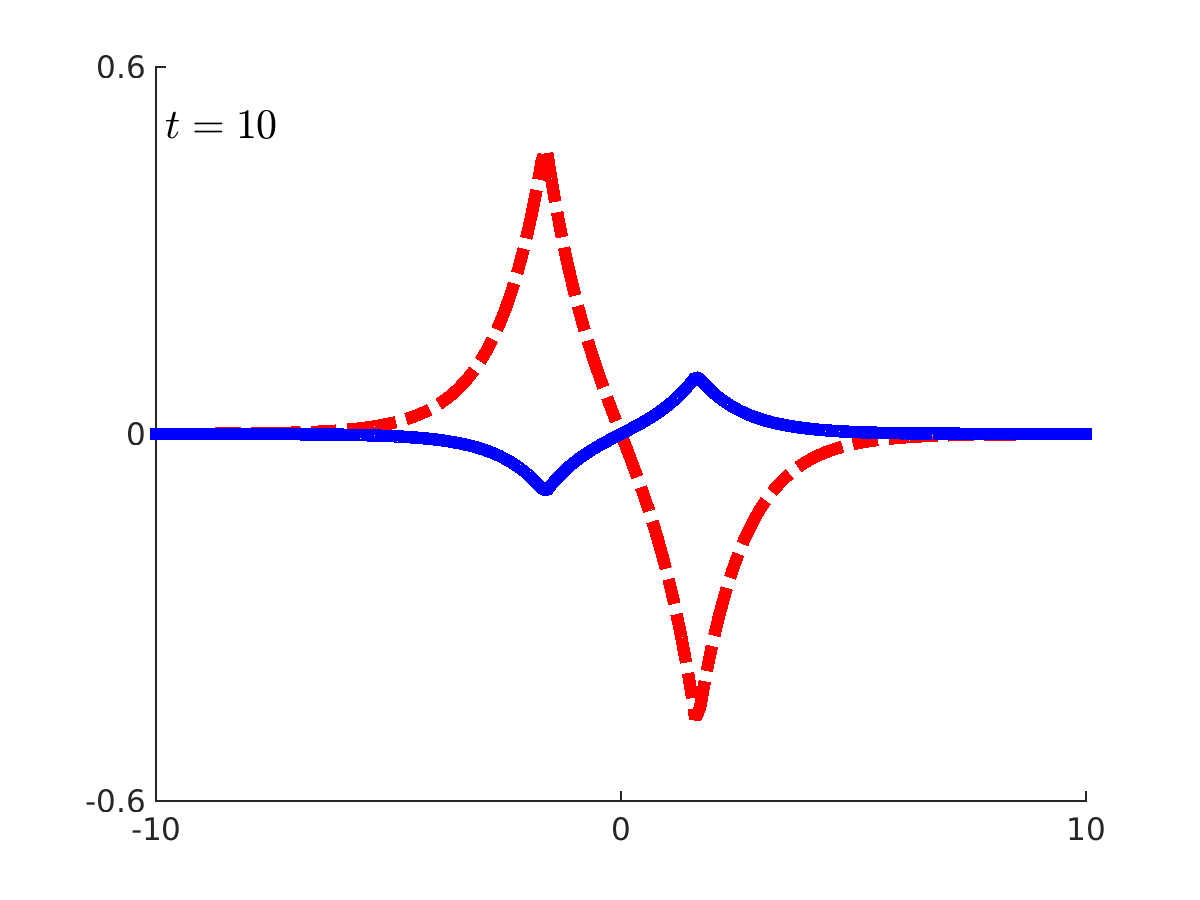}&\insertPlot{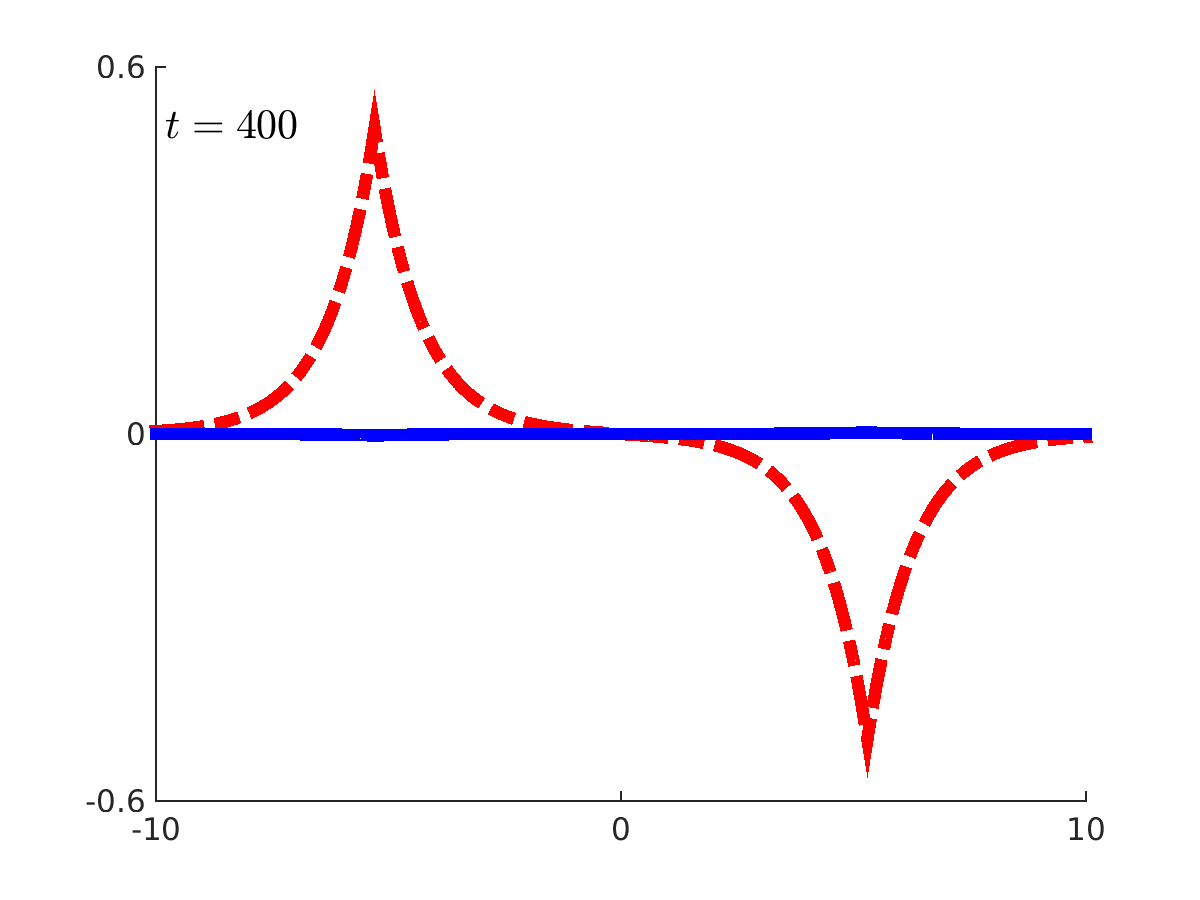}
  \end{tabular}
  \caption{Case $s = 1$. Limiting case. The solution decays to zero.}
  \label{fig:limcase}
\end{figure}

\section*{Acknowledgements}
K.G. gratefully acknowledges the hospitality of the {\em Institut Mittag--Leffler}, creating a great working environment for research during the fall 2016.


\begin{thebibliography}{10}

\bibitem{arnold1999topological}
V.~I. Arnold and B.A. Khesin.
\newblock {\em Topological methods in hydrodynamics}, volume 125 of {\em
  Applied Mathematical Sciences}.
\newblock Springer-Verlag, New York, 1998.

\bibitem{BressanFonte2005}
A.~Bressan and M.~Fonte.
\newblock An optimal transportation metric for solutions of the
  {C}amassa-{H}olm equation.
\newblock {\em Methods Appl. Anal.}, 12(2):191--219, 2005.

\bibitem{camassaholm1993}
R.~Camassa and D.~D. Holm.
\newblock An integrable shallow water equation with peaked solitons.
\newblock {\em Phys. Rev. Lett.}, 71(11):1661--1664, 1993.

\bibitem{constantin2000}
A.~Constantin.
\newblock Existence of permanent and breaking waves for a shallow water
  equation: a geometric approach.
\newblock {\em Ann. Inst. Fourier (Grenoble)}, 50(2):321--362, 2000.

\bibitem{constantinescher1998a}
A.~Constantin and J.~Escher.
\newblock Global existence and blow-up for a shallow water equation.
\newblock {\em Ann. Scuola Norm. Sup. Pisa Cl. Sci. (4)}, 26(2):303--328, 1998.

\bibitem{constantinescher1998b}
A.~Constantin and J.~Escher.
\newblock Wave breaking for nonlinear nonlocal shallow water equations.
\newblock {\em Acta Math.}, 181(2):229--243, 1998.

\bibitem{constantinescher2000}
A.~Constantin and J.~Escher.
\newblock On the blow-up rate and the blow-up set of breaking waves for a
  shallow water equation.
\newblock {\em Math. Z.}, 233(1):75--91, 2000.

\bibitem{constantin2008integrable}
A.~Constantin and R.~I. Ivanov.
\newblock On an integrable two-component {C}amassa-{H}olm shallow water system.
\newblock {\em Phys. Lett. A}, 372(48):7129--7132, 2008.

\bibitem{constantin2002geometric}
A.~Constantin and B.~Kolev.
\newblock On the geometric approach to the motion of inertial mechanical
  systems.
\newblock {\em J. Phys. A}, 35(32):R51--R79, 2002.

\bibitem{constantinstrauss2000}
A.~Constantin and W.~A. Strauss.
\newblock Stability of peakons.
\newblock {\em Comm. Pure Appl. Math.}, 53(5):603--610, 2000.

\bibitem{dikamolinet20092}
K.~El~Dika and L.~Molinet.
\newblock Stability of multi antipeakon-peakons profile.
\newblock {\em Discrete Contin. Dyn. Syst. Ser. B}, 12(3):561--577, 2009.

\bibitem{dikamolinet2009}
K.~El~Dika and L.~Molinet.
\newblock Stability of multipeakons.
\newblock {\em Ann. Inst. H. Poincar\'e Anal. Non Lin\'eaire},
  26(4):1517--1532, 2009.

\bibitem{grunert2015break}
K.~Grunert.
\newblock Blow-up for the two-component {C}amassa-{H}olm system.
\newblock {\em Discrete Contin. Dyn. Syst.}, 35(5):2041--2051, 2015.

\bibitem{grunert2012global}
K.~Grunert, H.~Holden, and X.~Raynaud.
\newblock Global solutions for the two-component {C}amassa-{H}olm system.
\newblock {\em Comm. Partial Differential Equations}, 37(12):2245--2271, 2012.

\bibitem{grunert2016}
K.~Grunert, H.~Holden, and X.~Raynaud.
\newblock Regularisation of the {H}unter--{S}axton and {C}amassa--{H}olm
  equations.
\newblock {\em (Submitted)}, 2016.

\bibitem{guantanwei2015}
C.~Guan, K.~Yan, and X.~Wei.
\newblock {Lipschitz metric for the modified two-component Camassa-Holm
  system}.
\newblock {arXiv:1510.03946}.

\bibitem{guan2010global}
C.~Guan and Z.~Yin.
\newblock Global existence and blow-up phenomena for an integrable
  two-component {C}amassa-{H}olm shallow water system.
\newblock {\em J. Differential Equations}, 248(8):2003--2014, 2010.

\bibitem{HoldenRaynaud2006}
H.~Holden and X.~Raynaud.
\newblock Convergence of a finite difference scheme for the {C}amassa-{H}olm
  equation.
\newblock {\em SIAM J. Numer. Anal.}, 44(4):1655--1680 (electronic), 2006.

\bibitem{HoldenRaynaud2006b}
H.~Holden and X.~Raynaud.
\newblock A convergent numerical scheme for the {C}amassa-{H}olm equation based
  on multipeakons.
\newblock {\em Discrete Contin. Dyn. Syst.}, 14(3):505--523, 2006.

\bibitem{holdenraynaud2007b}
H.~Holden and X.~Raynaud.
\newblock Global conservative multipeakon solutions of the {C}amassa-{H}olm
  equation.
\newblock {\em J. Hyperbolic Differ. Equ.}, 4(1):39--64, 2007.

\bibitem{holdenraynaud2007}
H.~Holden and X.~Raynaud.
\newblock Global conservative solutions of the {C}amassa-{H}olm equation---a
  {L}agrangian point of view.
\newblock {\em Comm. Partial Differential Equations}, 32(10-12):1511--1549,
  2007.

\bibitem{HoldenRaynaud2008p}
H.~Holden and X.~Raynaud.
\newblock A numerical scheme based on multipeakons for conservative solutions
  of the {C}amassa-{H}olm equation.
\newblock In {\em Hyperbolic problems: theory, numerics, applications}, pages
  873--881. Springer, Berlin, 2008.

\bibitem{holm2009singular}
D.~D. Holm, L.~{\'O}~N{\'a}raigh, and C.~Tronci.
\newblock Singular solutions of a modified two-component {C}amassa-{H}olm
  equation.
\newblock {\em Phys. Rev. E (3)}, 79(1):016601, 13, 2009.

\bibitem{kolev2004}
B.~Kolev.
\newblock Lie groups and mechanics: an introduction.
\newblock {\em J. Nonlinear Math. Phys.}, 11(4):480--498, 2004.

\bibitem{tan2011global}
W.~Tan and Z.~Yin.
\newblock Global periodic conservative solutions of a periodic modified
  two-component {C}amassa-{H}olm equation.
\newblock {\em J. Funct. Anal.}, 261(5):1204--1226, 2011.

\end{thebibliography}
\end{document}